# Mise en œuvre d'une ingénierie didactique de développement dans le cadre d'un travail collaboratif chercheur/enseignant lors de la conceptualisation des objets de l'Analyse au début du cursus dans le supérieur

## Fatma Belhaj Amor

**Abstract.** At the start of the higher education curriculum, the conceptualization of local approximation objects of a function requires the articulation of knowledge and skills from Functional Analysis and Topology. In the study of functions, a number of studies have established the existence of difficulties encountered by students, mainly as a result of the change of didactic contract during the transition from secondary to higher education. The construction of a teaching-learning project, as part of a collaborative effort with the class teacher, a priori helps students to overcome the main difficulties inherent in conceptualizing the local approximation objects of a function in the first year of preparatory classes. In the case of the design and implementation of didactic development engineering, analysis of the reasoning produced by students confronted with a situation with an adidactic dimension will enable us a priori to study the nature and origin of these difficulties. Our methodology for analyzing student work is based on a model of reasoning analysis within the framework of the theory of didactic situations mathematics. This model has played an essential role in the development of didactic engineering, in the identification of students' conceptions, forms and functions of reasoning. It also enabled us to identify epistemological, didactic and cultural obstacles to learning the concept of local approximation of a function.  These obstacles result either from the paradigm shift that takes place during the transition from secondary to higher education, or from working within the paradigm of Infinitesimal Analysis during the appropriation of this mathematical concept.

**Résumé.** Au début du cursus dans le supérieur, la conceptualisation des objets d'approximation locale d'une fonction nécessite l'articulation des connaissances et des savoirs qui relèvent de l'Analyse fonctionnelle et de la Topologie. Dans le domaine de l'étude des fonctions, certains travaux ont établi l'existence des difficultés rencontrées par les étudiants qui résultent principalement du changement de contrat didactique lors de la transition Secondaire/Supérieur. La construction d'un projet d'enseignement-apprentissage, dans le cadre d'un travail collaboratif avec l'enseignant de la classe, aide *a priori* les étudiants à surmonter les principales difficultés inhérentes à la conceptualisation des objets d'approximation locale d'une fonction en première année des classes préparatoires. Dans le cas de l'élaboration et la mise en œuvre d'une ingénierie didactique de développement, l'analyse des raisonnements produits par des étudiants confrontés à une situation à dimension adidactique nous permettra *a priori* d'étudier la nature et l'origine de ces difficultés. Notre méthodologie d'analyse des travaux des étudiants est basée sur un modèle d'analyse des raisonnements dans le cadre de la théorie des situations didactiques. Ce modèle a joué un rôle essentiel dans l'élaboration de l'ingénierie didactique de développement, dans l'identification des conceptions des étudiants, des formes et des fonctions des raisonnements. Il a permis aussi d'identifier des obstacles de nature épistémologique, didactique et culturelle liés à l'apprentissage du concept d'approximation locale d'une fonction. Ces obstacles résultent soit du changement de paradigme lors de la transition Secondaire/Supérieur, soit d'un travail s'inscrivant dans le paradigme *Analyse Infinitésimale* lors de l'appropriation de ce concept mathématique.

Table des matières

# 1. Introduction

Au début du cycle des classes préparatoires aux études d'ingénieurs, les objets d'approximation locale d'une fonction sont introduites afin d'être des nouvelles techniques permettant l'étude d'une fonction, le calcul de sa limite et l'existence de ses intégrales lors de la transition Secondaire/Supérieur, ainsi que l'étude de la convergence des séries dans le domaine des mathématiques en première année Physique-Chimie (PC). Par ailleurs, les étudiants mobilisent ces objets en tant que nouvelle méthode lors de la détermination des approximations de fonctions par des polynômes au voisinage d'un réel dans le domaine de la physique, l'optique, la mécanique, etc.

La conceptualisation des notions de fonctions négligeables, de fonctions équivalentes, de la formule de Taylor-Young et de développements limités nécessitent l'articulation des connaissances et des savoirs qui relèvent des domaines de l'Analyse fonctionnelle et de la Topologie. Par ailleurs, la notion de « voisinage d'un réel » est introduite dans le chapitre « Analyse asymptotique » de l'Analyse avant d'être traitée en Topologie à la fin de première année des classes préparatoires.

Dans le domaine de l'étude des fonctions, certains travaux ont établi l'existence des difficultés rencontrées par les étudiants qui résultent principalement du changement de contrat didactique lors de la transition Secondaire/Supérieur, (Pepin, 2014 ; Gueudet & Pepin, 2018 ; Bloch & Gibel, 2022). Ce changement de contrat se réalise :

- Soit au niveau du processus de l'étude des signes et des représentations (Bloch & Gibel, 2019). Dès le début de l'Université, les étudiants apprennent un nouveau signe comme on apprend un mot nouveau, selon Berger (2004). La culture symbolique attendue au niveau de l'Université diffère de celle développée au niveau Secondaire. Par exemple, un même symbole peut avoir de différents sens, (Bardini & Pierce, 2018) ;

- Soit au niveau de paradigme par l'existence des variables en transition liées au passage de « *l'Analyse concrète* » à « *l'Analyse moderne* » : obstacles épistémologiques et contrats, (Winslow, 2007) ;

-   Soit au niveau de la proposition des exercices aux élèves (Gueudet & Vandebrouck, 2019 ; Grenier-Boley, 2009) par le passage d'un « *travail algorithmique à un travail complexe* » (Ghedamsi, 2008). Et notamment, au niveau Secondaire, les tâches proposées ne laissent que peu d'autonomie aux élèves dans le choix des raisonnements (Leithner, 2000 ; Battie, 2003) ;

-   Soit au niveau de types des techniques de résolution des différentes tâches (Ghedamsi, 2008), le passage d'une « Analyse algébrisée » à une « Analyse formalisée » amène à l'existence des micro-ruptures d'ordres conceptuel et technique, (Praslon, 2000). L'approche algébrique fait obstacle lors de l'étude d'une fonction (Vandebrouck, 2011), et notamment lors de passage du niveau global au niveau local et réciproquement (Rogalski, 2008).

Dans ce cadre, le passage d'une "Analyse algébrisée" à une "Analyse formelle" peut amener à l'existence de « conflits conceptuels », de « micro-ruptures » ou d'obstacles de différentes natures. Des études soulignent la possibilité que certaines conceptions des étudiants entrant à l'Université constituent un obstacle à l'apprentissage des objets liés au champ de l'étude d'une fonction.

Par ailleurs, certains chercheurs explorent la manière adoptée par des mathématiciens pour transmettre certaines notions et certains concepts, comme la convergence ponctuelle dans une série de Taylor (Martin, 2012) et les notions de Topologie (Nemirovs & Smith, 2011 ; Wilkerson-Jerde & Wilensky, 2011).

Certains travaux conduits sur la notion de séries de Taylor ont établis que les étudiants sont confrontés aux difficultés pour réaliser les tâches (Alcock & Simpson, 2004, 2005) et des grandes difficultés à raisonner qui peuvent atteindre "un état de surcharge cognitive", (Kung & Speer, 2010).

Ces précédents travaux n'ont pas ciblé précisément l'analyse didactique des difficultés des étudiants, selon les aspects épistémologiques et cognitifs, lors de la conceptualisation des objets d'approximations locales des fonctions qui nécessite un travail complexe articulant les savoirs et connaissances qui relèvent des domaines de l'Analyse fonctionnelle et de la Topologie à l'entrée dans le Supérieur.

Ces constats nous amènent à nous interroger sur la nature et l'origine des difficultés rencontrées par les étudiants lors de l'appropriation et l'usage raisonné de ces objets mathématiques en classes préparatoires aux études d'ingénieures tunisiennes (IPEI), dans la section Physique-Chimie (PC).

Par ailleurs, certains travaux conduits sur la notion des séries de Taylor ont établi que les étudiants sont confrontés à des difficultés pour réaliser les tâches, (Alcock & Simpson, 2004, 2005), et des grandes difficultés à raisonner qui peuvent atteindre « un état de surcharge cognitive », (Kung & Speer, 2010). En effet,

> "Throughout the interviews, students took a very mechanical view of testing series, at times working through the steps of a particular test without being able to clearly articulate the underlying ideas, what they were testing for, or at times, even what the conclusion of the test was." (Kung & Speer, 2010, p. 9)

Kung et Speer (2010) interprètent les travaux d'Alcock en écrivant :

"…reported very little about students' responses to this question, in part because many students were simply confused about what the question was asking." (Ibid., p. 3)

Ces travaux nous amènent à penser qu'une expérimentation permettant le changement innovant du processus de l'enseignement du concept d'approximation locale d'une fonction aide les étudiants à surmonter les principales difficultés liées à l'appropriation et l'usage raisonné de ce concept mathématique en première année Physique-Chimie (PC). Ces difficultés résultent *a priori* en grande partie de l'absence des situations mathématiques dévolues aux étudiants permettant l'introduction des notions de fonctions négligeables, fonctions équivalentes, développements limités et la formule de Taylor-Young afin d'accéder au sens du savoir.

Dans le contexte de la confrontation des étudiants aux situations en classes ordinaires, le travail avec des enseignants extérieurs à la recherche, avant de la mise en œuvre d'une ingénierie didactique, permettra de mieux connaître ses pratiques ordinaires afin de surmonter des contraintes de l'enseignement ordinaire, (Perrin-Glorian & Baltar, 2019). Pour ces raisons, nous envisageons de construire et de mettre en œuvre, en étroite collaboration avec l'enseignante de la classe, une ingénierie didactique de développement en vue d'aider les étudiants à accéder à la définition formelle du concept d'approximation locale d'une fonction en première année (PC).

Dans le cadre des recherches portant sur l'interprétation des erreurs des élèves, Brousseau (1989) a prouvé, d'un côté, l'importance de la prise en compte des erreurs prévues dans la révolution historique des mathématiques pour son enseignement en précisant le travail du chercheur en didactique des mathématiques. D'un autre côté, il est très important dans les travaux de recherche d'un didacticien de connaître les difficultés des mathématiciens afin de les confronter à celles des étudiants. En effet,

"…il faut maintenant envisager les erreurs récurrentes comme le résultat (produit par et construit autour) de conceptions, qui, mêmes lorsqu'elles sont fausses, ne sont pas des accidents, mais des acquisitions souvent positives." (Brousseau, 1989, p.42)

A l'issue de ces travaux de Brousseau, nous envisageons d'élaborer et de construire une ingénierie didactique de développement qui nécessite de conduire préalablement une analyse épistémologique des concepts pour percevoir quel type d'obstacles peuvent apparaître et en même temps de questionner des obstacles qui peuvent être de nature didactique.

Dans le cadre de l'apprentissage des objets mathématiques, nous rejoignons le point de vue de Brousseau et Gibel (2005) qui donnent l'intérêt de l'analyse détaillée des raisonnements produits par les étudiants confrontés à une situation en vue d'étudier la nature et l'origine des erreurs apparaissant dans leurs productions.

A l'instar de Brousseau (1989), qui réfère aux travaux de Bachelard (1938) et Piaget (1975), les erreurs des élèves dans le domaine de la didactique peuvent relever en concept d'obstacle.

"Fondamentalement cognitifs, les obstacles semblent pouvoir être *ontogéniques*, *épistémologiques*, *didactiques* et même *culturels* selon leur origine et la façon dont ils évoluent." (Brousseau, 1989, p.3)

Et plus précisément, la notion d'obstacle épistémologique

"(…) tend à se substituer dans certains cas à celle d'erreur d'enseignement, d'insuffisance du sujet ou de difficulté intrinsèque des connaissances." (Brousseau, 1998, p.7-8)

Dans le contexte de l'analyse des difficultés rencontrées par les étudiants, nous rejoignons le point de vue de Brousseau (1986) qui donne l'intérêt d'une étude approfondie des différents types d'erreurs effectuées par les élèves confrontés à une situation mathématique. Cette analyse des erreurs se réalise par l'étude du développement cognitif et du fonctionnement des connaissances mobilisées par les élèves. Brousseau (1986) distingue deux catégories d'erreurs. D'une part, une erreur liée à une *connaissance fausse* qui peut être rejetée et surmontée comme son existence est passagère. Dans ce cas, elle reste une difficulté. D'autre part, un ensemble des connaissances fausses qui constitue une *conception erronée*. Cette conception sera difficile à modifier ou à rejeter et qui donnera lieu à des conflits cognitifs véritables. Dans ce cas, elle fait un obstacle : ontogénique, épistémologique lié à la genèse historique du concept ou didactique lié au choix de l'enseignement.

Nous proposons la question suivante :

Est-ce que dans le cas de l'élaboration d'une ingénierie didactique de développement, l'analyse des raisonnements produits par des étudiants confrontés à des situations réelles en classe permet d'identifier leurs conceptions (valides ou erronées) afin d'étudier la nature et l'origine de leurs difficultés inhérentes à la conceptualisation des objets d'approximation locale d'une fonction en première année (PC) ?

Dans ce travail, nous souhaitons que l'élaboration et la mise en œuvre d'une ingénierie didactique de développement (Perrin-Glorian, 2011), dans le cadre d'un travail collaboratif avec l'enseignante de la classe, lors de l'enseignement du chapitre « Analyse asymptotique », aide les étudiants à la conceptualisation des objets d'approximation locale d'une fonction.

Nous commençons cet article par la présentation de cadre et des outils théoriques adoptés afin de conduire une analyse des raisonnements produits par des étudiants (âgés de 18 à 20 ans) en termes de connaissances et savoirs mobilisés. Dans la deuxième partie, nous expliciterons la méthodologie de recherche mise en œuvre en présentant les différentes étapes successives de l'ingénierie. Ensuite, nous présenterons la méthode d'analyse de la séquence en mettant la focale sur la situation 1. La quatrième partie de cet article porte sur la présentation synthétique des principaux résultats de l'analyse *a posteriori* des travaux des étudiants qui est basée sur le modèle d'analyse du raisonnement introduit par Bloch et Gibel (2011). Finalement, nous présenterons l'intérêt de ce travail de recherche et nous ouvrirons autres perspectives dans le domaine de didactique de l'Analyse.

## 2. Cadre et outils théoriques

A l'instar de Gibel (2018), un raisonnement mathématique est classifié d'après ses fonctions et le type de situation étudiée. Il est intéressant de définir les conditions dans lesquelles les raisonnements mathématiques ont été produits lors de la résolution d'une situation adidactique ou à dimension adidactique au niveau du Primaire, du Secondaire et du Supérieur, (Ibidem). L'étude d'un raisonnement produit par un étudiant confronté à une situation d'action par la mise en considération du fonctionnement des connaissances nous permettra *a priori* d'identifier la nature et l'origine de ses difficultés.

### 2.A. La théorie des situations didactiques : cadre pertinent de l'analyse des raisonnements

Dans notre travail, nous mettons la focale sur l'étude du fonctionnement et du développement cognitif des connaissances et des savoirs mobilisés par l'étudiant lors de son action pour construire un apprentissage du concept d'approximation locale d'une fonction par confrontation aux situations mathématiques en classe. En mobilisant les concepts de la théorie des situations didactiques, l'analyse de la confrontation des étudiants à une situation mathématique nous offre la possibilité d'étudier leurs connaissances et leurs savoirs mobilisés, ainsi que le rôle de l'enseignant. Cette dialectique sujet/situation nécessite une étude des rapports établis entre l'étudiant et les savoirs inhérents à la Topologie et à l'Analyse fonctionnelle. L'analyse didactique de cette dialectique nous permet d'étudier les connaissances et les savoirs mobilisés (valides ou erronés) par l'étudiant en vue de produire un raisonnement.

Par ailleurs, le professeur a un rôle déterminant dans le processus enseignement-apprentissage des objets d'approximation locale d'une fonction en première année (PC). Il convient alors d'effectuer une analyse précise des interactions entre le sujet et le milieu avec lequel il interagit. Nous devons conduire une analyse didactique de la relation entre enseignant, étudiants, savoirs et milieu effectuée lors de la confrontation des étudiants à une situation réelle en classe. Cette relation didactique se réalise dans un champ de l'étude des expérimentations. La théorie des situations didactiques (Brousseau, 1998) est un cadre théorique pertinent qui permet l'analyse de ces expérimentations par l'adaptation des notions de situation, d'institutionnalisation des savoirs et des connaissances et schéma de la structuration du milieu (Bloch, 2005).

### 2.A.a. Notion de situation à dimension adidactique

Dans le système sujet-milieu, l'activité de l'enseignant et de l'étudiant s'articulent tout au long de déroulement d'une séquence. En effet, en classe lors de l'interaction de l'étudiant en situation réelle, il peut se saisir des problèmes et des difficultés, à ce moment, l'enseignant doit interagir pour lui aider. Et plus précisément,

> "Du point de vue de l'activité mathématique de l'enseignant et de l'élève, il apparaît bien qu'elles sont conjointes, ou plutôt dépendantes l'une de l'autre : elles s'articulent durant tout le déroulement de la séance, elles s'alimentent l'une l'autre." (Bloch, 1999, p.33)

Dans notre cas de la confrontation des étudiants à une situation à dimension adidactique liée aux notions d'approximation locale d'une fonction, l'enseignant fournit à la réflexion de ses étudiants. Dans ce cas, on parle, (Bloch, 1999), d'une situation à dimension adidactique qui

"- d'une part, est dépendante des connaissances présentes dans le milieu et des connaissances que l'enseignant peut y injecter ;

- d'autre part se manifeste par une interaction de connaissances élèves / professeur."
(Bloch, 1999, p.33)

Dans le cas d'une situation à dimension adidactique, l'enseignant admet un rôle important lors des interactions étudiant/savoirs et étudiant/milieux.

Dans notre travail, nous nous intéressons notamment à une situation d'action au point de vue de la construction d'un savoir.

Une séquence d'enseignement intégrant une ou plusieurs situations à dimension adidactique est fréquemment organisée en trois « phases » :

- phase d'action : la connaissance du sujet se manifeste par des décisions : succès ou échec ;

- phase de formulation : elle met en rapport deux sujets avec le milieu, et le succès exige que l'un des sujets formule la connaissance à l'intention de l'autre ;

- phase de validation théorique : construction formelle dans un répertoire de règles ou de théorèmes reconnus.

### 2.A.b. Institutionnalisation des savoirs et connaissances en termes de répertoire didactique

Pour Brousseau (2010), l'institutionnalisation des savoirs est liée à un processus didactique et résulte d'une intervention spécifique permettant à l'étudiant et l'enseignant de « reconnaître » et de légitimer « l'objet d'enseignement ».

Lorsque l'étudiant est confronté à une situation d'action, il a recours à ses connaissances antérieures et ses savoirs qui sont sous forme de définitions, propriétés, formules, théorèmes, des règles du calcul formel et algébrique, etc. Cet ensemble des savoirs et des connaissances mobilisés constitue le répertoire didactique, (Gibel, 2018).

- <u>Répertoire didactique de la classe</u>

Dans le cas de l'enseignement de la notion d'approximation locale d'une fonction en première année (PC), le répertoire didactique de la classe est constitué des attendus des étudiants par le professeur, suite à son enseignement du concept d'approximation locale d'une fonction dans le chapitre « Analyse asymptotique ». Dans le champ de l'étude locale d'une fonction, nous pouvons distinguer deux composantes du répertoire didactique :

- le répertoire didactique du Secondaire : il intègre toutes les propriétés, les formules, les techniques, les théorèmes, les règles du calcul algébrique et formel enseignés au niveau Secondaire ;

- le répertoire didactique du Supérieur : il intègre les propriétés, les formules, les définitions, les techniques, les théorèmes, etc. qui sont liés aux notions de fonctions équivalentes, de fonctions négligeables, de voisinage d'un réel, de développement limité d'une fonction au voisinage d'un point et à la formule de Taylor-Young. En fait, il est le répertoire didactique des connaissances et des savoirs enseignés dans la classe de première année (PC).

- Répertoire didactique de l'étudiant

En classes préparatoires, les étudiants considèrent que les cours et les travaux dirigés dispensés par l'enseignant ne suffisent pas pour l'apprentissage et l'usage raisonné des notions et des concepts mathématiques. Nous rejoignons le point de vue de Farah (2015) et Lalaude (2016) qui indiquent que les étudiants ont fréquemment recours, en dehors du temps scolaire, à l'utilisation de certains ouvrages universitaires pour investiguer de nouvelles situations, inhérentes à la notion ou au concept visé en vue d'approfondir leurs connaissances et leurs savoirs. D'une certaine manière, ceci constitue pour les étudiants un prolongement de l'étude. En conséquence, l'étudiant, confronté à une situation mathématique en classe, a recours à son propre répertoire didactique, selon Gibel (2018).

La résolution par l'étudiant des situations adidactiques et à dimensions adidactique nécessite l'utilisation de son répertoire de représentation. De ce fait, chaque étudiant décide de la mise en œuvre d'une suite d'actions, valides ou erronées, relevant de l'usage de son répertoire d'action. À partir des résultats obtenus, il va modifier son répertoire de représentation, (Gibel, 2018). Cette interaction entre les répertoires est résumée dans un schéma permettant de modéliser le fonctionnement du répertoire didactique.

"Dans la modélisation du fonctionnement des connaissances, en théorie des situations (…), celles-ci apparaissent comme les moyens hypothétiques, pour le sujet, de prendre des décisions." (Gibel, 2018, p.31)

Nous considérons le schéma théorique, proposé par Gibel et Ennassef (2012), du fonctionnement des connaissances de l'étudiant confronté à une situation adidactique (d'action) où à dimension adidactique.

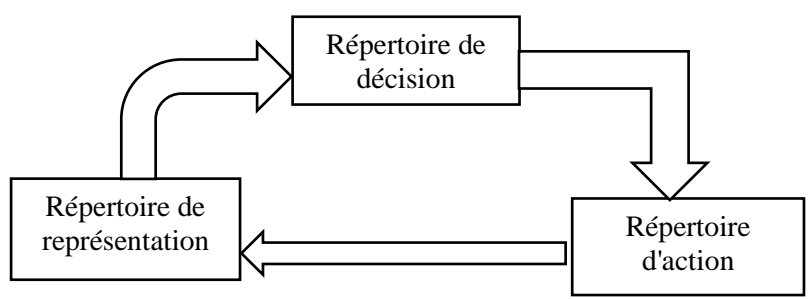

Figure 1 - Schéma de modélisation du fonctionnement du répertoire didactique de l'étudiant en situation d'apprentissage (Gibel & Ennassef, 2012)

Dans notre travail, ce schéma nous permettra de modéliser le fonctionnement de connaissances de l'étudiant en situation réelle en classe. Cette modélisation est réalisée par la précision des composantes du répertoire didactique de la classe et du répertoire didactique de chaque étudiant (des théorèmes, formules, propriétés, symboles, techniques utilisées, etc.), selon Gibel (2018).

Pour ceci, nous adoptons, dans le cadre de la théorie des situations didactiques, le modèle d'analyse des raisonnements introduit par Bloch et Gibel (2011) et nous envisageons de l'enrichir par des outils théoriques permettant l'étude des connaissances et des savoirs d'ordres topologique et fonctionnel.

## 2.B.  Outillages didactiques pour analyser la construction d'un concept de l'Analyse

### 2.B.a. Caractérisation d'une conception mathématique

Dans le contexte de l'analyse des erreurs des étudiants, à l'instar de Brousseau (1998), le sens d'une connaissance mathématique se définit par plusieurs moyens parmi eux l'ensemble des conceptions.

> "Le sens d'une connaissance mathématique se définit, non seulement par la collection des situations où cette connaissance est réalisée en tant que théorie mathématique, (…), non seulement par la collection des situations où le sujet l'a rencontrée comme moyen de solution, mais aussi par l'ensemble des conceptions, des choix antérieurs qu'elle rejette, des erreurs qu'elle évite, les économies qu'elle procure, les formulations qu'elle reprend, etc." (Brousseau, 1998, p.117)

Dans notre travail, une conception est un outil d'ordre théorique qui permet l'analyse des observations des étudiants en situation, (Artigue, 1991).

> "… mais l'identification de conceptions locales qui se manifeste en situation et l'analyse des conditions de passage de telle conception locale à telle autre, qu'il s'agisse de rejeter une conception erronée, de mettre en place une conception permettant d'améliorer l'efficacité dans la résolution de telle ou telle classe de problèmes ou de favoriser la mobilité entre des conceptions déjà disponibles. De ce point de vue, c'est l'objet local qui est bien l'outil adéquat." (Artigue, 1991, p.271-272)

De ce fait, une *conception* est considérée comme un outil pour la construction d'un concept mathématique qui permet la modélisation de l'étudiant en tant qu'agent apprenant, (Ibidem). Par ailleurs, une *conception* est liée à la connaissance de l'apprenant en situation. Elle est distinguée par les représentations et les traitements qu'elle mobilise dans un domaine scolaire. En conclusion, une « *conception est l'instance d'un concept par le système du [sujet<>milieu]* », (Ibid.).

En s'attachant à la nature de la situation, Brousseau a toujours voulu indiquer qu'on ne pouvait pas étudier les connaissances et les savoirs mobilisés sans tenir compte des conditions dans lesquelles ils sont produits. C'est-à-dire qu'il faut toujours mettre au centre de la situation les questions suivantes : Pourquoi pose-t-on ce type de problème ? Quelle est la nature de ce type de problème ?

Quand on a bien analysé la nature de la situation, on peut montrer les conceptions qui sont susceptibles de surgir et les conceptions qui surgirent effectivement. En adoptant le paradigme de l'erreur et la problématique des obstacles, Balacheff (1995) considère qu'une conception

> "… n'a pas en soit de caractère local, elle est d'abord une connaissance au sens d'un état d'équilibre dynamique du système en interaction sujet/milieu." (Balacheff, 1995, p.224)

Le système des représentations joue un rôle déterminant dans la conception. Nous allons attacher à définir les représentations selon la sémiotique de Peirce et les registres de représentation sémiotique de Duval (1993, 2006).

## 2.B.b. Représentations sémiotiques d'un concept mathématique

Nous pensons que l'étude des signes nous aidera à conduire une analyse plus fine d'un concept mathématique. La *sémiotique* désigne l'étude des signes et de leurs significations.

La définition de l'étude des signes « l'analyse sémiotique » est introduite par Peirce au XIX^e siècle et développée par Duval en fin du XX^e siècle. Nous commençons cette partie par la présentation de la sémiotique peircienne. La théorie de Peirce est explicitée par Bloch (2005). En effet,

" "Un signe a, comme tel, trois références : premièrement il est un signe pour quelque pensée qui l'interprète ; deuxièmement, il est un signe qui tient lieu de quelque objet auquel il est équivalent dans cette pensée ; troisièmement il est un signe, sous quelque rapport ou qualité, qui le met en connexion avec cet objet. Posons-nous la question de savoir ce que sont les trois corrélats auxquels une pensée-signe réfère". (1868)

"Un signe est quelque chose qui représente une autre chose pour un esprit. Pour son existence comme tel trois choses sont requises. En premier lieu, il doit avoir des caractères qui nous permettront de le distinguer des autres objets. En second lieu, il doit être affecté d'une façon ou d'une autre par l'objet qui est signifié ... En troisième lieu chaque représentation s'adresse à un esprit. C'est seulement pour autant qu'elle fait ceci qu'elle est une représentation" (1873) (…)." (Bloch, 2005, p.9)

Par ailleurs, une analyse de l'activité mathématique nécessite des outils permettant l'analyse des connaissances et des savoirs mobilisés par l'étudiant. Comme indiqué dans les travaux de Bloch et Gibel (2011), cette analyse nécessite d'avoir recours à une analyse sémiotique de Peirce afin d'analyser les signes produits, leurs usages et leurs transformations. A l'instar de ces didacticiens, l'étude des signes qui apparaissent dans la production des étudiants nous aide à identifier les différentes conceptions, valides ou erronées.

Comme le spécifie Duval (2006), la distinction entre un objet et ses différentes représentations est obligatoire dans l'analyse sémiotique d'une activité mathématique.

"En mathématiques, une représentation n'est intéressante que dans la mesure où elle peut se transformer en une autre représentation. Un signe n'est intéressant que s'il peut être substitué à d'autres signes pour effectuer des opérations (Condillac, 1982)." (Duval, 2006, p.57)

Le seul moyen d'avoir accès à un objet mathématique est de prendre en compte ses différentes représentations sémiotiques qui

"… sont des productions constituées par l'emploi de signes appartenant à un système (sémiotique) de représentation qui a ses propres contraintes de signifiance et de fonctionnement. Une figure géométrique, un énoncé en langue naturelle, une formule algébrique, un graphe sont des représentations sémiotiques qui relèvent de systèmes sémiotiques différents." (Duval, 1993, p.39)

Comme l'indique Hitt (2004), il est important de prendre en considération la production des représentations sémiotiques par les étudiants dans une démarche heuristique.

"Le caractère fonctionnel et dynamique de ces représentations sémiotiques permet à l'étudiant d'avoir un contrôle sur la démarche qui le mène à la solution." (Hitt, 2004, p.350)

Dans notre étude des objets inhérents à l'approximation locale d'une fonction, nous pouvons distinguer plusieurs registres de représentation sémiotique :

- Le registre algébrique : il intègre des expressions algébriques (polynômes, équation de la tangente, etc.), de fonctions algébriques, des transformations des formules algébriques, etc.

- Le registre analytique : il concerne les expressions algébriques qui comprennent des symboles (exemples : $\varepsilon, \Sigma$ , o, etc.), des expressions ou représentations analytiques comme les concepts de l'Analyse (limite, dérivée, intégrale, continuité, prolongement par continuité, etc.). Par exemple, le développement limité d'une fonction au voisinage de 0 à l'ordre n s'écrit sous la forme d'une égalité : $f(x) = P(x) + x^n \varepsilon(x)$ Avec $\lim_{x \to 0} \varepsilon(x) = 0$ et $P$ : un polynôme de $IR_n[X]$ .

- Le registre graphique : il comporte les courbes représentatives des fonctions (algébriques, transcendantes, etc.), les tangentes en un point, etc.

- Le registre numérique : il comprend tout ce qui est calcul numérique comme les approximations décimales, rationnelles, calcul numérique de limite d'une fonction, à l'aide de la calculatrice, d'une fonction au voisinage ou en un point, calcul des dérivées en un point, etc.

- Le registre géométrique : il désigne toute interprétation des illustrations graphiques ou des expressions analytique ou algébrique comme l'étude locale de la position relative de la courbe par rapport à sa tangente ou à son asymptote, etc.

Comme indiqué dans les travaux de Bloch et Gibel (2011), l'analyse des raisonnements nécessite d'effectuer une analyse sémiotique afin d'étudier les signes produits, leurs usages et leurs transformations. Ainsi, la sémiotique triadique de Peirce, couplée à l'approche sémiotique de Duval (2006), pourrait nous permettre d'analyser le processus des raisonnements produits par les étudiants.

## 2.B.c. Evolution des objets de l'Analyse standard à la transition Secondaire/Supérieur – Les différents paradigmes liés à l'apprentissage de l'Analyse

L'étude historico-épistémologique du concept d'approximation locale d'une fonction, nous a permis d'établi l'existence des obstacles (Brousseau, 1989) qui résultent d'une rupture d'ordre épistémologique et sont étroitement liés au changement de paradigme (Belhaj Amor, 2022). Dans ce contexte, nous rejoignons le point de vue de Kuzniak, Montoya, Vandebrouck et Vivier (2015) qui ont conduit une étude de l'évolution historique des mathématiques et son impact sur l'enseignement actuel en identifiant trois paradigmes de l'Analyse standard :

"Paradigme [Analyse Arithmetico-géométrique] (AG) (…) qui permet des interprétations provenant, avec quelques implicites, de la géométrie, du calcul arithmétique mais aussi du monde réel (…).

Paradigme [Analyse Calculatoire] (AC) : (…) dans ce calcul algébrique généralisé, les règles de calcul sont définies, plus ou moins explicitement, et elles sont appliquées indépendamment d'une réflexion sur l'existence et la nature des objets introduits (…).

Paradigme [Analyse Infinitésimale] (AI) : "(…) un travail spécifique et formel s'appuie sur l'approximation et la localité : bornes, inégalités, travail sur des voisinages, négligeabilité, (…)." (Kuzniak, Montoya, Vandebrouck & Vivier, 2015, p.7)

Dans la pratique de l'Analyse en première année (PC), nous pouvons envisager dans le cas de l'étude locale d'une fonction au voisinage d'un réel :

$\exists \propto > 0, \forall x \in ] -\alpha; \alpha[ \cap D_f$ ; on a : $f(x) = a_0 + a_1 x + a_2 x^2 + o(x^2)$ un développement limité d'ordre 2 de la fonction $f$ au voisinage de 0.

- *Paradigme [AG]* : par exemple, le réel $a_1$ est la pente de la tangente à la courbe représentative de la fonction $f$ au point $(0, f(0))$ qui est d'équation : $y = a_0 + a_1 x$.

- *Paradigme [AC]* : par exemple, le résultat d'un calcul symbolique sur des expressions analytiques, le calcul de limite du rapport et de différence des fonctions, etc.

- *Paradigme [AI]* : par exemple, la fonction $g(x) = a_0 + a_1 x + a_2 x^2$ est l'approximation polynômiale d'ordre 2 de la fonction $f$ au voisinage de 0.

Dans notre étude, nous nous intéressons à identifier clairement dans quel paradigme s'inscrit le raisonnement de l'étudiant et notamment sa manière d'appréhender les conceptions. Nous pouvons distinguer quatre catégories de conceptions :

- *conceptions spécifiques* qui relèvent du paradigme [AC] et/ou [AG] : dans ce cas, l'étudiant effectue un calcul algébrique et formel ou il mobilise des objets de géométrie. Par exemple, il détermine l'équation d'une tangente par la mobilisation de la technique du calcul du nombre dérivé ou par l'interprétation de sa représentation graphique, etc.

- *conceptions plus fines et plus précises* qui relèvent du paradigme [AI] : dans ce cas, l'étudiant mobilise les notions de fonctions négligeables, d'approximations polynômiales d'une fonction, de développements limités ou de la formule de Taylor-Young afin de déterminer l'approximation locale d'une fonction ;

- *conceptions* s'inscrivent dans le paradigme [AI] qui coexistent à un certain moment de la progression d'un raisonnement en lien avec les paradigmes [AC] et [AG] : dans ce cas, l'étudiant s'appuie sur des objets du domaine de l'étude d'une fonction en utilisant le concept d'approximation locale d'une fonction. Par exemple, il interprète le résultat du calcul de limite de fonctions par la mobilisation de la technique de la détermination de fonctions négligeables ;

- *conceptions* que certains étudiants ne maîtrisent pas suffisamment : dans ce cas, leurs conceptions ne reflètent pas la mobilisation de connaissances et de savoirs spécifiques du Supérieur qui s'inscrivent dans le paradigme [AI].

Nous pensons que l'identification des paradigmes de l'Analyse standard nous permettra d'étudier d'une manière plus précise l'évolution des connaissances et des savoirs mobilisés par les étudiants confrontés à une situation réelle en classe en vue de produire un raisonnement.

Lorsque l'étudiant va être confronté à une situation à dimension adidactique, nous allons intéresser au point de vue de conceptions de connaître si sa production s'inscrivant dans quel paradigme de l'Analyse.

## 2.B.d. Analyse logique d'un concept de l'Analyse en termes de dimension sémantique et/ou syntaxique

Afin de conduire une analyse plus approfondie des raisonnements, nous pensons nécessaire de relier le signe, son sens et sa dénotation lors de la construction et l'utilisation du concept d'approximation locale d'une fonction. Pour cela, nous procéderons à une analyse logique d'un concept mathématique en termes de dimensions sémantique et/ou syntaxique introduites dans les travaux de Frege (1971) et Tarski (1974) et développées dans les travaux (Kouki, 2018 ; Kouki, Belhaj Amor & Hachaïchi, 2016 ; Belhaj Amor, 2016, 2018, 2019, 2020) pour étudier les raisonnements produits par les étudiants confrontés à une situation à dimension adidactique.

Dans le cas de l'analyse logique de notre objet d'étude, nous distinguons la dimension sémantique, la dimension syntaxique et l'articulation de ces deux dimensions lors de la construction des conceptions valides dans des productions des étudiants. La syntaxe fournit des règles des transformations et formulations des expressions analytique et algébrique lors de son appropriation et de son usage par l'adaptation des définitions ou théorèmes ou propriétés du concept mathématique. La sémantique désigne toute interprétation des représentations graphiques des objets d'approximation locale d'une fonction, ainsi que toute interprétation ou vérification des expressions algébriques ou analytiques. L'articulation sémantique/syntaxique est liée à la mobilisation des moyens de contrôle, de justification et de vérification lors de l'utilisation et la mise en œuvre d'une propriété ou d'une formule, d'une définition ou d'un théorème, etc. du répertoire didactique.

## 2.B.e. Niveaux de mise en fonctionnement des connaissances par les étudiants

Au début du cursus dans le Supérieur, lorsque les étudiants sont confrontés en classe à des situations mathématiques, nous conduisons des analyses didactiques de leurs productions en mettant en jeu l'usage d'une connaissance ou d'un savoir. Dans ce contexte, nous rejoignons le point de vue de Robert (1998) qui a déterminé les caractéristiques de mise en fonctionnement des connaissances et des savoirs, en termes d'outil, relatives à un niveau scolaire bien déterminé. Elle a distingué trois niveaux de mises en fonctionnement des connaissances par les élèves :

"Le niveau technique : (…) Ce niveau correspond pour nous à des mises en fonctionnement indiquées, isolées, mettant en jeu des applications immédiates de théorèmes, propriétés, définitions, formules, etc. " (Robert, 1998, p.165)

"Le niveau des connaissances mobilisables : ce niveau correspond à des mises en fonctionnement plus large : encore indiquée, mais dépassant l'application simple d'une propriété à la fois. (…)" (Ibid., p.166)

"Le niveau des connaissances disponibles : ce niveau correspond au fait de savoir résoudre ce qui est proposé sans indications, d'aller rechercher soi-même dans ses connaissances ce qui peut intervenir." (Ibid., p.167)

Ce dernier niveau est plus spécifique *a priori* au Supérieur. Il nous permet de connaître la manière de l'évolution des pratiques attendues des étudiants lors de la production d'un raisonnement.

## 3. Méthodologie générale de recherche

La méthodologie mise en œuvre repose sur l'élaboration d'une ingénierie didactique de développement dans le cadre d'un travail collaboratif avec l'enseignante de la classe pendant deux années universitaires (2018-2019 et 2019-2020). Cette méthodologie de recherche se réalise en trois étapes successives que nous présenterons d'une manière synthétique dans cette section (Belhaj Amor, 2022).

### 3.A. Etape 1 : phase préliminaire

Une étude historico-épistémologique a permis de conclure que la genèse historique du concept d'approximation locale d'une fonction est due à l'usage des représentations graphiques et géométriques. De plus, il existe des obstacles épistémologiques qui résultent de changement de paradigmes de "cinématique" au géométrique puis analytique.

Une étude institutionnelle du programme officiel de première année (PC) et de polycopié de cours de l'enseignante nous amène à conclure que les enseignants sont invités à introduire une nouvelle notion à partir de la mobilisation de ses différentes représentations (algébrique, analytique, graphique et géométrique) par la prise en considération des connaissances antérieures de l'étudiant. Par contre, l'approche graphique n'est pas convoquée par l'enseignante, (Ibidem). Elle suive l'enseignement "classique" lors de l'introduction des définitions formelles des objets d'approximation locale d'une fonction dans le chapitre "Analyse asymptotique".

L'analyse des raisonnements produits par des étudiants confrontés à une situation mathématique « classique », dans le cadre d'une évaluation formative et formatrice, en avril 2019, nous a permis d'identifier et de caractériser des obstacles épistémologique et didactique liés à l'apprentissage des objets d'approximation locale d'une fonction en première année (PC), (Belhaj Amor & Gibel, 2022). Ces obstacles résultent *a priori* en grande partie de l'absence des situations mathématiques dévolues aux étudiants permettant l'introduction du concept d'approximation locale d'une fonction qui nécessite l'articulation des éléments d'ordre topologique et fonctionnel dans le paradigme [AI].

### 3.B. Etape 2 : la construction du chapitre « Analyse asymptotique »

Nous construisons et élaborons une ingénierie didactique de développement en collaboration étroite avec l'enseignante de la classe qui a un rôle très spécifique dans ce projet d'enseignement-apprentissage du concept d'approximation locale d'une fonction. Elle co-construit le chapitre

« Analyse asymptotique » en intégrant deux situations à dimension adidactique en vue d'offrir aux étudiants la possibilité d'accéder au sens du savoir "développement limité d'une fonction au voisinage d'un réel" par un travail s'inscrivant dans les différents paradigmes de l'Analyse standard. La construction de ce concept mathématique se réalise à partir du fonctionnement de leurs connaissances antérieures, par leur action sur le milieu.

### 3.C.  Etape 3 : élaboration du modèle d'analyse du raisonnement

Une fois cette ingénierie est mise en œuvre, notre méthodologie des analyses *a priori* et *a posteriori* des travaux des étudiants est basée sur le modèle d'analyse des raisonnements introduit par Bloch et Gibel (2011) et développé dans notre travail de recherche. Nous envisageons de l'enrichir par des outils théoriques permettant une analyse plus fine des raisonnements, au point de vue de conceptions, en vue d'identifier les connaissances et les savoirs mobilisés dans quel paradigme de l'Analyse standard (Kuzniak, Montoya, Vandebrouck & Vivier, 2015).

Dans notre étude, nous conduirons une analyse des raisonnements qui est orientée selon quatre axes dont le premier est lié au niveau du milieu de la situation et de paradigmes correspondants. Le deuxième axe se réalise selon l'analyse des fonctions des raisonnements (Bloch & Gibel, 2011). Le troisième est l'étude des signes (Peirce) et des représentations sémiotiques (Duval, 1993) observables. Le dernier axe est l'étude de l'usage logique des connaissances antérieures de l'étudiant et les éléments du contrôle. Nous présentons la synthèse de la méthodologie d'analyse des raisonnements dans ce schéma :

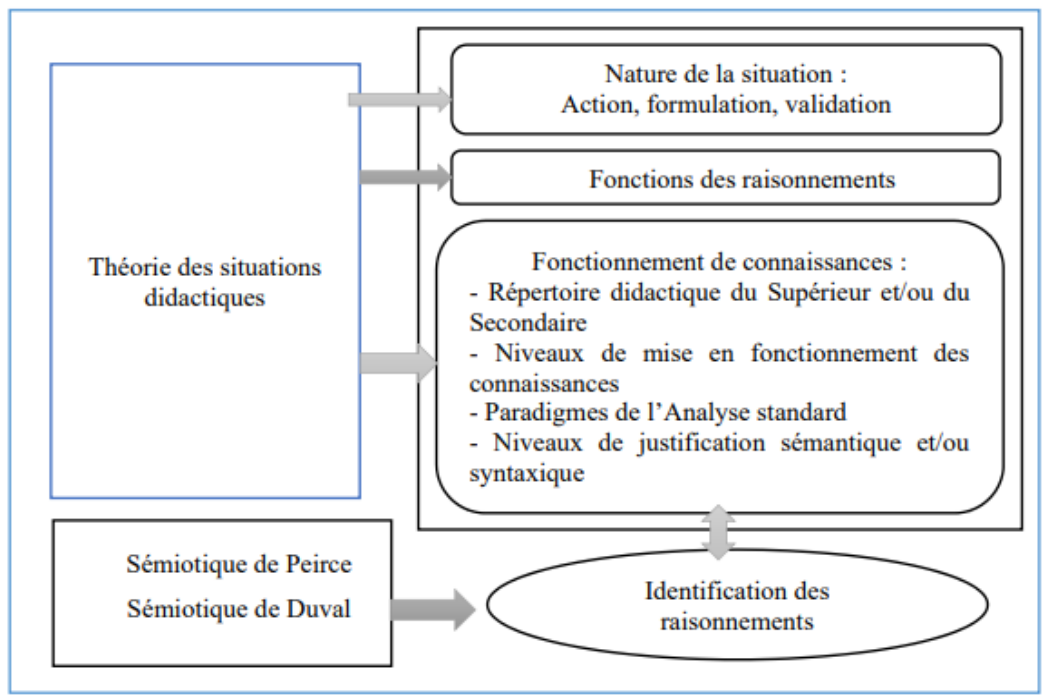

Figure 2 - Schéma du modèle d'analyse des raisonnements

Dans notre travail, nous adoptons le cadre de la théorie des situations didactiques en mettant la focale sur les notions de conception, connaissance, savoir, situation mathématique, situation à dimension adidactique, institutionnalisation et structuration du milieu. Lorsque les étudiants vont se confronter à une situation à dimension adidactique, nous allons nous intéresser à l'analyse de leur production en termes de raisonnements :

- du point de vue de conceptions : c'est de connaître si la production relève du paradigme [AC], si elle relève du paradigme [AG], si elle relève du type d'Analyse spécifique *a priori* au niveau Supérieur : le paradigme [AI] ;

- du point de vue logique : c'est un travail spécifique en deux éléments très importants du raisonnement par l'identification de la dimension sémantique ou de la dimension syntaxique ou de l'articulation des deux dimensions ;

- du point de vue des signes : c'est l'étude des signes qu'apparaissent dans un raisonnement par l'analyse sémiotique de Peirce en référant aux registres de représentation sémiotique de Duval (1993) ;

- du point de vue du fonctionnement des connaissances et des savoirs mobilisés par leurs actions sur le milieu : c'est l'étude de fonction, forme et nature du raisonnement.

## 4. La séquence objet d'étude : focus sur la situation 1

### 4.A. Sujets

Dans l'étude réalisée, les étudiants sont inscrits en première année des classes préparatoires de la filière Physique-Chimie (PC) et ont été admis au baccalauréat Sciences Expérimentales avec un score élevé.

Une enseignante expérimentée[1] de l'Institut Préparatoire aux Études d'Ingénieurs de Tunis (IPEIT), madame Rajia Slim, a accepté de participer à ce travail de recherche. Cette enseignante a montré un vif intérêt à s'appuyer sur les résultats de l'évaluation de ses étudiants en vue d'améliorer le contenu de ses enseignements.

### 4.B. Instrumentations

#### 4.B.a. Les données recueillies avant la mise en œuvre de la situation

L'enseignante a planifié son enseignement du chapitre « Analyse asymptotique » de manière à programmer deux séances afin de proposer une situation en classe. Elle a accepté d'interagir avec nous et de suivre nos recommandations inhérentes à la gestion des phases au cours desquelles les étudiants sont confrontés à une situation à dimension adidactique en privilégiant, pour répondre à leurs éventuelles questions, le recours à une guidance faible. Elle organise le processus de la transposition didactique des savoirs à enseigner aux savoirs enseignés de ce chapitre en suivant le déroulement prévu de notre expérimentation.

Par ailleurs, elle a déposé une demande d'autorisation de filmer les séances et de nous permettre de participer à l'expérimentation en classe. L'accord du responsable administratif de l'Institut Préparatoire (IPEIT) nous permet d'enregistrer les séances en des vidéos et de prendre des photos.

---

[1] Ayant enseigné depuis de plus de dix ans dans cette filière. Elle est en poste depuis sa réussite au concours du recrutement dans l'enseignement Supérieur des docteurs de mathématiques.

L'enseignante commence le chapitre « Analyse asymptotique » par l'introduction du concept de la relation de comparaison des suites et des fonctions. A la fin de cette séance, elle informe ses étudiants des deux classes de (PC) qu'ils vont travailler en groupes hétérogènes, plus précisément en trinômes, pendant les deux prochaines séances dans l'amphi. Ils ont accepté que nous récupérions la photocopie de leurs productions après que l'enseignante les ait informés des enjeux et du cadre de la recherche. Notre expérimentation se déroule pendant la deuxième et la troisième séance de l'enseignement de ce chapitre. La situation 1 est proposée, en janvier 2020, en classe de première année (PC) comportant 24 étudiants qui est divisée en groupes de trois étudiants.

Notre expérimentation se déroule pendant la deuxième et la troisième séance de l'enseignement de ce chapitre.

### 4.B.b. Les données recueillies lors de la mise en œuvre de la situation

L'enseignante commence la séance par un rappel sur les notions de fonctions dominées, de fonctions négligeables, de voisinage d'un réel et de fonctions équivalentes. Ensuite, elle distribue elle-même à chaque étudiant l'énoncé de la situation 1 afin de leur laisser d'avantage du temps pour travailler en trinôme. Tous les documents sont autorisés. Quand les étudiants rencontrent des difficultés pour répondre aux questions et que leurs échanges ne leur permettent pas de les surmonter alors ils ont la possibilité de demander l'aide auprès de l'enseignante pour pallier les insuffisances du milieu auquel ils sont confrontés. Cette dernière, par une guidance faible (Bartolini-Bussi, 2009), le plus souvent sous la forme d'un questionnement adéquat, pourra les amener à envisager une interprétation de leurs productions ou à prendre de nouvelles décisions. Après quarante minutes, l'enseignante ramasse les copies des trinômes. Par la suite, un représentant de chaque trinôme vient au tableau pour rédiger la réponse du trinôme à l'une des questions posées. L'intervention de l'enseignante dépend de la nature de la réponse, valide ou erronée, rédigée par le représentant du trinôme. Son rôle principal est de susciter le débat au sein de la classe quant à la pertinence, la validité et la complexité du raisonnement produit. Ensuite, elle va s'appuyer sur ces raisonnements produits par les étudiants afin d'introduire la formule de Taylor-Young.

### 4.C. Présentation des enjeux de la situation

Dans notre étude, nous construisons une situation à dimension adidactique dont l'objectif principal consiste à amener les étudiants à comprendre qu'à une fonction $f$ quelconque (algébrique non polynômiale ou transcendante), de classe $C^n$ au voisinage d'un réel $x_0$, on peut déterminer à différents degrés des approximations polynômiales et locales au voisinage de ce réel. La structuration de cette situation mathématique est similaire à celle étudiée dans l'ouvrage de terminale (1983), (Belhaj Amor, 2022). D'une part, nous voulons prendre en considération les connaissances antérieures des étudiants en permettant de déterminer eux-mêmes les approximations locales d'ordre 1, 2 et 3 d'une fonction via ses approximations polynômiales successives de degré 1, 2 et 3 au voisinage d'un réel. D'autre part, nous souhaitons leur permettre d'introduire eux-mêmes ce nouveau concept mathématique en articulant les dimensions sémantique et syntaxique par un travail s'inscrivant dans les différents paradigmes de l'Analyse.

- <u>Énoncé de la situation mathématique</u>

On considère les fonctions suivantes :

$f(x) = \sqrt{1+x}, f_1(x) = 1 + \frac{x}{2}, \quad f_2(x) = 1 + \frac{x}{2} - \frac{x^2}{8}$ et $f_3(x) = 1 + \frac{x}{2} - \frac{x^2}{8} + \frac{x^3}{16}$

**1.a)** Calculer $\lim_{x \to 0}[f(x) - f_1(x)]$

 **b)** Calculer $\lim_{x \to 0} \frac{f(x) - f_1(x)}{x}$

 **c)** En déduire l'expression de la fonction $f$ en fonction de $f_1$ au voisinage de 0.

**2.a)** Calculé $\lim_{x \to 0} \frac{f(x) - f_2(x)}{x^2}$

 **b)** En déduire l'expression de $f$ en fonction de $f_2$ au voisinage de 0.

**3.a)** Calculer $\lim_{x \to 0} \frac{f(x) - f_1(x)}{x^2}$.

 **b)** En déduire l'expression de $f$ en fonction de $f_1$ au voisinage de 0.

(En déduire $f(x) = f_1(x) - \frac{1}{8}x^2 + o(x^2)$)

 **c)** Préciser la position de la courbe $C_{f_1}$ par rapport à celle de $C_f$ au point $(0; f(0))$.

**4)** Montrer que $f_2$ est la seule fonction polynôme de degré 2 tel que : $\lim_{x \to 0} \frac{f(x) - f_2(x)}{x^2} = 0$

**5.a)** D'une manière analogue à la question 2, déterminer l'expression de $f$ en fonction de $f_3$ au voisinage de 0.

 **b)** Préciser la position de la courbe $C_{f_2}$ par rapport à celle de $C_f$ au point $(0; f(0))$.

**6.a)** Calculer $f(0), f'(0), f''(0)$ et $f'''(0)$.

 **b)** Exprimer $f_1$ en fonction de $f(0)$ et $f'(0)$.

 **c)** Exprimer $f_2$ en fonction de $f(0), f'(0)$ et $f''(0)$.

 **d)** En déduire l'expression de $f_3$ en fonction de $f(0), f'(0), f''(0)$ et $f'''(0)$.

## 4.D. Eléments d'analyse *a priori* de la situation

Dans cette partie, nous présentons des éléments d'analyse *a priori* de la situation sur le plan didactique, ainsi que le déroulement prévu de la séquence.

### 4.D.a. Analyse *a priori* de la situation sur le plan didactique

La situation est assimilable à une situation à dimension adidactique en vue d'introduire les approximations locales d'une fonction d'ordre 1, 2 et 3 à partir d'un cas particulier afin de formaliser le concept de la formule de Taylor-Young.

Les objectifs de la proposition de cette situation consistent à mobiliser et réorganiser des connaissances antérieures de l'étudiant, ainsi qu'à introduire les approximations d'une fonction en faisant le lien entre ces connaissances et le nouvel objet mathématique.

Nous pouvons distinguer trois variables didactiques :

**VD$_1$** : La nature de la fonction étudiée et ses approximations polynômiales au voisinage de 0.

**VD$_2$** : Le degré des approximations polynômiales étudiées.

**VD$_3$ :** Le point au voisinage duquel les étudiants doivent produire l'approximation polynômiale locale (dont le degré est fixé).

### 4.D.b. Analyse en termes de niveaux de milieux et de paradigmes de l'Analyse

Cette situation est une situation d'activité permettant d'introduire la formule de Taylor-Young. Nous rejoignons le point de vue de Bloch et Gibel (2011) qui ont distingué les niveaux de milieux de $M_0$ à $M_{-2}$ dans le cas de l'étude d'une situation adidactique ou à dimension adidactique.

"… car c'est au niveau de l'articulation entre le milieu objectif (milieu heuristique) et le milieu de référence que nous pourrons voir apparaitre et se développer les raisonnements attendus." (Bloch & Gibel, 2011, p.10)

Les questions successives visent à introduire des approximations polynômiales locales d'une fonction à l'ordre 1, puis à l'ordre 2 et enfin à l'ordre 3 au voisinage de 0 à partir du calcul de limites de fonctions développé dans le paradigme [AC] et à percevoir les conditions définissant les fonctions négligeables dans le paradigme [AI]. Nous essayons de voir les échanges au sein du groupe qui amènent à un niveau de confrontation en situation de référence. Ainsi, nous sommes en train de poser des questions permettant aux étudiants de passer de la situation de référence à la situation d'apprentissage, ainsi que de changement de paradigmes de [AC] à [AI], puis de [AI] à [AG].

- **$S_{-2}$ : Situation de référence : sujet agissant et milieu objectif**

Le milieu objectif est constitué des expressions analytique et algébrique, ainsi que des types de tâches à effectuer afin de déterminer l'approximation polynômiale de la fonction $f$ au voisinage de 0 successivement à l'ordre 1, puis à l'ordre 2, puis à l'ordre 3. L'étudiant va établir une action sur des objets mathématiques qui sera motivée par son répertoire didactique. En effet, il va mobiliser ses connaissances antérieures du répertoire didactique pour le calcul de limites des fonctions dans le paradigme [AC]. Nous rejoignons le point de vue de Gibel et Ennassef (2010) : l'élève va élaborer un répertoire de décision de stratégie du calcul de limite de la fonction définie par un rapport des fonctions.

"Ainsi, par l'usage de leur répertoire de représentation, ils décideront de la mise en œuvre d'une suite d'actions sur le milieu matériel. Cette suite d'actions, valides ou erronées, relève de leur répertoire didactique, plus précisément du fonctionnement de leur système organisateur." (Gibel & Ennassef, 2010, p.3)

Les actions de l'étudiant agissant sur le milieu objectif ont pour objet le calcul de limites des rapports de deux fonctions, ainsi que la manière de l'effectuer, de l'obtenir et de le contrôler dans le paradigme [AC].

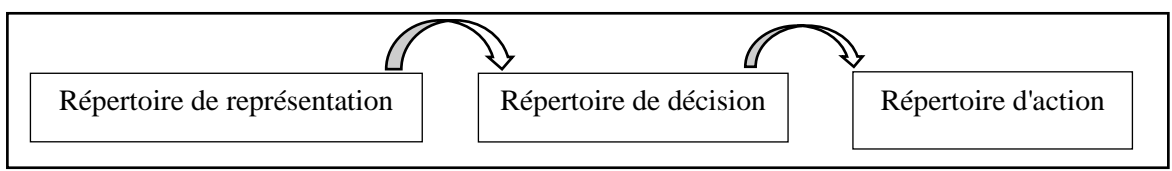

Figure 1 - Modélisation du fonctionnement du répertoire didactique de l'élève (Gibel & Ennassef, 2010, p.4)

Les expressions analytiques et algébriques données vont alors prendre le statut d'argument au niveau de milieu objectif. Et plus précisément, les expressions $(\lim_{x \to 0} \frac{f(x) - f_1(x)}{x} = 0,\ \lim_{x \to 0} \frac{f(x) - f_1(x)}{x^2} = -\frac{1}{8},\ \lim_{x \to 0} \frac{f(x) - f_2(x)}{x^2} = 0\ et\ \lim_{x \to 0} \frac{f(x) - f_3(x)}{x^3} = 0)$ vont prendre les statuts d'arguments au niveau de milieu objectif en mettant en jeu les conditions nécessaires pour l'usage de fonctions négligeables.

- **$S_{-1}$ : Situation d'apprentissage : étudiant apprenant et milieu de référence**

L'étape qui suit le calcul de limite, réalisée par un travail s'inscrivant dans le paradigme [AC], correspond à la capacité de l'étudiant à interpréter ce résultat en référence aux fonctions négligeables par un usage raisonné de l'expression du reste : en $\varepsilon(x)$ ou $o(x^n)$.

On a : $\lim_{x \to 0} \frac{F(x)}{G(x)} = 0$ donc on en déduit en recours au cours que $F(x) = o(G(x))$ ou $\frac{F(x)}{G(x)} = \varepsilon(x)$ avec $\lim_{x \to 0} \varepsilon(x) = 0$.

La situation d'apprentissage repose sur la capacité de l'étudiant à :

- réaliser un lien avec le cours sur les fonctions négligeables et à en déduire les expressions analytiques attendues. En effet, il est capable d'exprimer la fonction $f$ en fonction de $f_1$ puis en fonction de $f_2$, puis en fonction de $f_3$, ainsi que de percevoir le comportement des fonctions $(f - f_1)$, $(f - f_2)$, $(f - f_3)$, etc. L'étudiant va utiliser son répertoire didactique, plus précisément, le système générateur en ayant pris conscience du fait que les conditions de son usage sont vérifiées ;

- prouver l'existence des approximations locales d'une fonction au voisinage de 0 d'ordre 1, 2 et 3 et son unicité ;

- comprendre que le passage d'une approximation polynômiale d'ordre 1 à une approximation polynômiale d'ordre 2 délivre des informations très intéressantes pour la détermination de la position relative de deux courbes représentatives de cette fonction et de son approximation polynômiale locale d'ordre 1 ;

- réaliser un lien entre la fonction $f$ définie au voisinage de 0, sa valeur en ce réel, les valeurs de ses dérivées première, seconde et troisième en ce réel et les approximations polynômiales proposées. Les étudiants doivent percevoir les conditions nécessaires pour réaliser le calcul des dérivées successives d'une fonction en 0.

Les étudiants vont s'interroger sur la validité de chacune des expressions analytiques trouvées dans les paradigmes [AC] puis [AI], ainsi que sur la position relative des deux courbes dans les paradigmes [AI] puis [AG].

L'enseignante a pour objectif de prendre en considération le travail produit par les étudiants en vue d'institutionnaliser la formule de Taylor-Young.

- **$S_0$ : Situation didactique : étudiant confronté au milieu $M_0$**

L'enjeu didactique de cette situation repose sur la manière de mettre en forme l'approximation d'une fonction au voisinage d'un point et la mise en généralisation. Et notamment, l'enjeu est de comprendre qu'une fonction n-fois dérivable au voisinage d'un point peut être approximée par un polynôme. En effet, nous pouvons « approcher », de plus en plus près, la fonction étudiée, autrement dit avec une précision croissante liée au degré du polynôme.

Nous allons expliciter une manière d'établir une approximation polynômiale locale d'une fonction de classe $C^n$ au voisinage d'un réel donné par la recherche des nombres dérivés successifs en ce réel.

La situation repose sur l'idée de comprendre la manière d'affiner progressivement l'approximation d'une fonction, au voisinage de 0, par des polynômes. Cette approximation se fait tout d'abord par une fonction polynôme de degré 1, puis par une fonction polynôme de degré 2 qui est obtenue en ajoutant un terme à la précédente, puis un travail sur la fonction polynôme de degré 3 qui est obtenue en ajoutant de nouveaux termes de degré 3. À partir de ces expressions, l'enseignante va ensuite introduire la formule de Taylor-Young dans le paradigme [AI].

**4.D.c. Déroulement prévu de la séquence**

L'enseignante commencera la séance par un rappel écrit détaillé sur le concept de la relation de comparaison de fonctions, et plus précisément sur les notions de fonctions négligeables, de fonctions dominées et de fonctions équivalentes en lien étroit avec la notion de voisinage. Ensuite, l'enseignante va dévoluer la situation 1. Elle donnera une feuille pour chaque étudiant contenant la situation en leur demandant de travailler en groupe pendant quarante minutes. Elle ramassera les copies des étudiants. Ensuite, un représentant de chaque groupe présentera, au tableau, leur production pour chacune des questions proposées.

Après la présentation par chaque représentant d'un groupe de leur raisonnement au tableau, l'enseignante interviendra et à ce moment-là il y aura une interaction et un débat en classe. Son intervention dépend de la réponse de ce groupe (valide ou erroné). L'enseignante s'appuiera sur la sixième question, l'unicité du développement limité d'une fonction au voisinage d'un réel et les différentes écritures syntaxiques des approximations polynômiales locales de la fonction $f$ pour donner la formule de Taylor-Young.

**4.E. Méthode d'analyse des données**

Notre méthodologie repose sur la prise en compte des quatre composantes du modèle d'analyse de raisonnement, développé précédemment, qui est déduit du modèle d'analyse des raisonnements de Bloch et Gibel (2011) et de Gibel (2018), lors de l'analyse des différents types de raisonnement élaborés par les étudiants dans une situation à dimension adidactique. Par ailleurs, dans ce type de situation les niveaux de milieux où se situent les étudiants sont : le milieu objectif ($M_{-2}$), le milieu de référence ($M_{-1}$) et le milieu d'apprentissage ($M_0$).

Cette méthodologie d'analyse du raisonnement par la mise en jeu de ses composantes (fonctions, dimension sémiotique, nature et justification) associée à chaque niveau de milieux $M_{-2}$, $M_{-1}$ et $M_0$ est réalisée selon ces quatre axes :

- Fonctions des raisonnements : dans chaque trinôme, les étudiants vont interagir entre eux pour décider des actions à réaliser sur les objets mathématiques, effectuer et contrôler le calcul de limites de fonctions, de différence de deux fonctions et les rapports de fonctions, réaliser et contrôler le calcul des valeurs des dérivées successives de la fonction $f$ en 0 et établir l'existence et l'unicité de la fonction $f_2$, l'approximation polynômiale locale d'ordre 2 de la fonction $f$ au voisinage de 0 dans le paradigme [AI].

- Niveau d'utilisation des symboles, des signes et des représentations sémiotiques : les étudiants vont devoir mobiliser des signes et des symboles pour effectuer des calculs de différences et de rapports des fonctions par un travail s'inscrivant dans le paradigme [AC].

- Utilisation et caractérisation du répertoire didactique selon les niveaux de paradigme : les étudiants vont mobiliser leurs connaissances antérieures relatives au système générateur et au registre des formules. Et plus précisément, ils vont utiliser des méthodes de calcul de limites de fonctions et le calcul des dérivées successives de la fonction $f$ dans le paradigme [AC], ainsi que l'étude de la position relative des courbes représentatives des deux fonctions dans le paradigme [AG].

    - Nature de la justification d'un raisonnement (les dimensions sémantique et/ou syntaxique) : les étudiants vont justifier leurs raisonnements en vérifiant les conditions nécessaires et suffisantes pour l'usage des signes et des différentes formules du répertoire didactique. Par exemple, lors de la démonstration de l'unicité d'approximation polynômiale locale d'ordre 2 de la fonction $f$ au voisinage de 0, ils vont faire l'hypothèse de l'existence d'une autre approximation polynômiale locale d'ordre 2 de cette fonction au voisinage de ce réel afin de déduire la contradiction avec les hypothèses de la situation mathématique.

Les étudiants ont produit des calculs et des raisonnements en ayant recours à leur répertoire didactique. Selon le contenu des questions de la situation, les actions de l'étudiant agissant ont pour objet de produire un pas du raisonnement de différentes natures :

- Pas de raisonnement de nature syntaxique : dans ce cas, l'étudiant va mobiliser des règles du calcul formel et algébrique, des propriétés, des formules et des théorèmes pour produire le pas de raisonnement objet d'étude.

- Pas de raisonnement de nature sémantique : dans ce cas, l'étudiant va produire un pas de raisonnement en justifiant et/ou en interprétant les expressions analytiques du calcul de limites des fonctions. Par exemple, il va interpréter les nouvelles représentations analytiques de la fonction $f$ et plus précisément, ses approximations locales d'ordre 2 et 3 au voisinage de 0 afin de préciser la position relative de deux courbes en un point d'abscisse 0.

- Pas de raisonnement articulant les dimensions sémantique et syntaxique : dans ce cas, l'étudiant va justifier et vérifier les conditions de l'utilisation d'une technique du répertoire didactique.

    Concernant les types de raisonnements produits par les étudiants, nous adoptons la caractérisation du point de vue de Gibel (2004) :

Raisonnement formel complet : dans ce cas, les étudiants donnent du sens à la connaissance mathématique mobilisée. Par exemple, ils formulent le développement limité d'une fonction, précise l'ordre de développement limité et le voisinage sur lequel il est défini.

Raisonnement formellement incomplet : dans ce cas, les étudiants produisent des procédures de résolutions qui ne mentionnent pas explicitement les conditions de l'utilisation des connaissances antérieures utilisées (c'est-à-dire la validité des hypothèses nécessaires lors de l'utilisation d'une formule ou d'un énoncé). D'une part, ils ne précisent pas l'intervalle de voisinage sur lequel une approximation locale de la fonction $f$ est définie. D'autre part, il ne vérifie pas que cette fonction est de classe $C^{\infty}$ au voisinage de 0 pour déterminer les valeurs de ses dérivées successives en ce réel.

Raisonnement erroné : dans ce cas, les étudiants produisent des procédures de résolutions fausses.

Absence de raisonnement : dans ce cas, les étudiants ne donnent aucune réponse.

L'objectif principal de cette analyse des raisonnements produits par les trinômes consiste à identifier leurs difficultés rencontrées lors de la conceptualisation du concept d'approximation locale d'une fonction au voisinage de 0 en première année (PC).

Nous conduisons une analyse des raisonnements produits par les 8 trinômes en distinguant les différents types de techniques mobilisées lors de la résolution des questions proposées. Nous mettons en jeu les différents paradigmes de l'Analyse (Kuzniak, Montoya, Vandebrouck & Vivier, 2015) afin d'identifier les techniques mobilisées par les étudiants lors de la résolution des différentes questions proposées.

A cause des contraintes éditoriales, dans cet article, nous choisissons de présenter les principaux résultats de ces analyses en adoptant le modèle d'analyse des raisonnements.

## 5. Principaux résultats de l'analyse des raisonnements

Après la distribution des feuilles contenant l'énoncé de la situation mathématique à toute la classe, les étudiants travaillent individuellement, puis ils décident de produire un écrit qui est la synthèse des échanges sur leurs procédures de résolutions. À certains moments, les étudiants au sein du groupe étaient en désaccord sur leurs productions, certains d'entre eux ont demandé l'aide auprès de l'enseignante. Dans ce cas, l'enseignante a essayé de les guider pour mobiliser les connaissances nécessaires pour produire un raisonnement valide. De ce fait, les étudiants ont produit un raisonnement « commun » qui résulte des interactions au sein du trinôme. Afin d'analyser les raisonnements produits par ces étudiants, nous nous appuierons sur l'analyse *a priori* et le déroulement prévu réalisés auparavant.

Dans une première étape, nous conduisons une analyse quantitative des productions de l'ensemble des trinômes pour chacune des questions posées dans la situation. Ensuite, grâce à l'enregistrement vidéo réalisé, nous mettons la focale sur l'analyse des productions de l'un des trinômes qui ont demandé l'aide auprès de l'enseignante afin d'identifier la nature et l'origine de leurs difficultés relatives à la conceptualisation des objets d'approximation locale d'une fonction. Dans ce cas, nous choisissons de conduire une analyse qualitative détaillée des raisonnements produits par le trinôme 1 avant et après l'intervention de l'enseignant, ainsi que de certains épisodes relatifs à la formulation de leur procédure. Cette analyse joue un rôle important afin d'étudier des difficultés liées à l'appropriation et l'usage raisonné du concept de fonctions négligeables.

### 5.A. Analyse globale des raisonnements

Dans le cadre de la théorie des situations didactiques, nous conduisons une analyse de notre corpus, constitué des productions de 8 trinômes et des retranscriptions des interactions en classe, en termes de raisonnements, selon le modèle d'analyse développé précédemment. Nous présentons les principaux résultats de cette analyse quantitative en les classant selon les objectifs de la proposition des différentes questions de la situation. À l'issue de cette identification, nous pouvons présenter une catégorisation des différents types de techniques utilisées en mettant en jeux les registres et les paradigmes associés au point de vue de conceptions :

- Technique Arithmetico-algébrique : par exemple, l'étudiant utilise la méthode de l'expression conjuguée pour calculer la limite du rapport de deux fonctions. Son travail s'inscrit dans le paradigme [AC].

- Technique Arithmetico-analytique : par exemple, il prouve que la fonction $f$ est de classe $C^3$ au voisinage de 0 afin de calculer les valeurs de ses dérivées successives en ce réel dans le paradigme [AC].

- Technique analytique : par exemple, l'étudiant utilise le théorème de l'Hôpital pour calculer la limite du rapport de deux fonctions dans le paradigme [AC]. Il fait appel aussi à la méthode de raisonnement par l'absurde pour prouver l'unicité d'approximation locale d'une fonction au voisinage de 0.

- Technique géométrique : par exemple, il étudie le signe de la différence de deux fonctions pour préciser la position relative de ses deux représentations graphiques dans le paradigme [AG].

- Technique analytique-infinitésimale : par exemple, il utilise la formule de fonctions négligeables du répertoire didactique du Supérieur afin de déterminer une approximation locale d'une fonction au voisinage de 0 par un travail s'inscrivant dans le paradigme [AI].

Nous présentons le résumé des différents types de de raisonnements et de techniques mobilisés par les 8 trinômes dans le tableau ci-dessous :

| | Types de techniques utilisées | Raisonnements formels complets | Raisonnements formellement incomplets | Raisonnements erronés | Absence de raisonnement |
|---|---|---|---|---|---|
| Calcul de limite de fonctions | Arithméco-algébrique | 4 | 0 | 3 | 0 |
| | Analytique | 0 | 1 | 0 | |
| Détermination d'approximations locales de $f$ | Analytique-infinitésimale | 1 | 6 | 0 | 0 |
| | Autres | 0 | 1 | 0 | 0 |
| Étude de la position de deux courbes | Géométrique | 3 | 0 | 2 | 2 |
| | Aucune | 0 | 0 | 1 | |
| Preuve de l'unicité d'approximation locale d'ordre 2 de $f$ au V(0) | Analytique | 1 | 0 | 0 | 7 |
| Calcul des dérivées de $f$ en 0 | Arithméco-analytique | 6 | 0 | 1 | 1 |
| Détermination d'approximations polynômiales de $f$ par l'usage des symboles de ses dérivées successives en 0 | Arithméco-analytique | 4 | 0 | 2 | 0 |
| | Aucune | 2 | 0 | 0 | |

Tableau 3 - Types de raisonnement et techniques mobilisés lors de la résolution des différentes questions proposées dans la situation mathématique

Nous présentons les principaux résultats de l'analyse de notre corpus selon les objectifs des questions proposées.

-   Calculus de limite de fonctions (questions (1-a-b), (2-a) et (3-a))

La plupart des trinômes ont produit des pas de raisonnements de nature syntaxique lors du calcul de limite de fonctions en utilisant une technique Arithmetico-algébrique par un travail s'inscrivant dans le paradigme [AC]. Par ailleurs, la plupart d'entre eux sont confrontés aux problèmes de justification sémantique des différents pas de leur raisonnement au niveau des connaissances mobilisables.

Nous remarquons que tous les étudiants ont mobilisé des connaissances et des savoirs du répertoire didactique du Secondaire lors du calcul de limite du rapport de deux fonctions dans le cas d'une forme indéterminée. Ils ne parviennent pas à mobiliser des connaissances du répertoire didactique du Supérieur et plus précisément le théorème de l'Hôpital. Nous pouvons dire qu'il existe une rupture conceptuelle qui s'opère lors l'enseignement des concepts de l'Analyse entre le niveau Secondaire et le niveau Supérieur lié au champ de limite. Ainsi, il existe *un obstacle de nature culturelle* par le recours aux connaissances vues au Secondaire lors de la résolution des questions traitées au lycée.

-   Détermination des approximations polynômiales locales d'ordre 1, 2 et 3 de la fonction $f$ au voisinage de 0 (questions (1-c), (2-b) et (3-b))

Concernant l'interprétation des résultats du calcul de limite, la majorité des trinômes ont produit des raisonnements formellement incomplets. En effet, ils ont déterminé les expressions analytiques des approximations polynômiales locales d'ordre 1, 2 et 3 de la fonction $f$ par un travail s'inscrivant en réalité dans le paradigme [AC]. En effet, ils ne précisent pas les voisinages de 0 pour lesquels doivent être déterminées ces approximations. Ils sont confrontés à une difficulté liée aux connaissances disponibles "fonctions négligeables". Ils ont utilisé une technique analytique-infinitésimale du répertoire didactique du Supérieur sans connaître le sens de la connaissance mathématique mobilisée « approximation locale d'une fonction ». De ce fait, les étudiants sont confrontés à une difficulté liée à la définition formelle du concept de la relation de comparaison de fonctions. En effet, l'écriture d'une fonction sous forme d'une somme d'un polynôme et d'un reste ne donne pas du sens à ce concept mathématique. Ces difficultés résultent, en grande partie, à la complexité de l'articulation des éléments d'ordre topologique et fonctionnel lors de la détermination des différentes écritures syntaxiques du concept de fonctions négligeables qui nécessitent des moyens du contrôle et de validation de son existence selon le point de vue sémantique. Ces erreurs commises apparaissent lors du passage de niveau du milieu objectif au niveau du milieu de référence, ainsi que lors du changement de paradigmes [AC] au paradigme [AI] par l'interprétation des résultats des calculs effectués dans le milieu $M_{-2}$. Nous pouvons dire que ces étudiants ont effectué des *erreurs d'ordre cognitif et conceptuel* lors de la conceptualisation de la notion de relation de comparaison de fonctions. Ces erreurs relèvent d'un *obstacle de nature didactique*.

Par ailleurs, le groupe 4 est le seul trinôme qui a produit des pas de raisonnements de nature sémantique/syntaxique, par une guidance faible, lors de la détermination des approximations polynômiales locales d'ordre 1, 2 et 3 de la fonction $f$ au voisinage de 0.

-   Etude de la position relative de deux courbes (questions (3-c) et (5-b))

Certains trinômes ont produit des pas de raisonnements de nature sémantique/syntaxique lors de l'usage d'une technique géométrique afin de préciser la position relative de deux courbes dans le paradigme [AG]. En effet, ils ont étudié le signe de la différence de deux fonctions $f$ et $f_1$ pour préciser la position relative de ses représentations graphiques en utilisant des connaissances mobilisables dans des registres analytique et géométrique.

Par contre, deux trinômes (groupe 1 et 4) considèrent le signe de la représentation analytique de la fonction $[x \mapsto f(x) - f_1(x)]$ dépend du signe de celle de $[x \mapsto \varepsilon(x) - \frac{1}{8}]$. L'étude du signe de cette différence prouve l'existence d'une difficulté liée à la compréhension de la définition de la fonction $[x \mapsto \varepsilon(x)]$ donnée lors de l'introduction du concept de relation de comparaison de fonctions (fonctions négligeables et fonctions équivalentes au voisinage d'un réel). Cette *erreur d'ordre cognitif et conceptuel* résulte de la dialectique technique/sens de ce concept mathématique lors de changement de paradigme de [AI] à [AG]. Nous pouvons dire que cette erreur relève d'un *obstacle lié aux choix didactiques de l'enseignement*. En effet, l'enseignante définie la fonction $[x \mapsto \varepsilon(x)$, dans son cours, lors de l'appropriation du concept de relation de comparaison de fonctions dans le paradigme [AI] juste par l'expression suivante : $[\lim_{x \to 0} \varepsilon(x) = 0]$.

Nous remarquons qu'aucun trinôme n'a construit les représentations graphiques des deux fonctions $f$ et $f_1$ afin de préciser leur position relative. Par contre, ils ont l'habitude d'interpréter graphiquement ses représentations, au niveau Secondaire, pour répondre à ce type de question.

  - Preuve de l'unicité d'approximation locale d'ordre 2 de la fonction au voisinage de 0 (question (4))

Un seul trinôme (groupe 4) a produit des pas de raisonnement de nature sémantique par l'usage de la méthode du raisonnement par l'absurde du répertoire didactique du Secondaire. Ce groupe a développé cette technique analytique dans le paradigme [AI] afin de produire un raisonnement formel complet. Par contre, les autres trinômes n'ont pas trouvé la méthode adéquate pour prouver l'unicité de l'approximation locale d'ordre 2 de la fonction au voisinage de 0. Nous pouvons dire qu'il existe une difficulté persistante depuis le Secondaire lors de la résolution de ce type de tâche qui est liée à la démonstration d'une preuve mathématique : l'existence et l'unicité d'un objet de l'Analyse.

  - Calcul des dérivées de la fonction (question (6-a))

La majorité des trinômes ont calculé les dérivées successives de la fonction $f$ sans apporter les moyens de la justification nécessaire de l'usage de la technique analytique du répertoire didactique du Supérieur, ce qui garantit l'existence de ses dérivées successives en ce réel. Par ailleurs, un seul groupe a produit des éléments d'ordre sémantique : la fonction $f$ est de classe $C^{\infty}$ sur ]-1,+∞[. Leurs travaux s'inscrivent dans le paradigme [AC] par la mobilisation d'une technique algébrique dans des     registres algébrique et numérique. L'usage automatique des formules du répertoire didactique de la classe prouve que les étudiants n'ont pas l'habitude de justifier leur raisonnement.

  - Détermination des approximations polynômiales de la fonction au voisinage de 0 par l'usage des symboles de ses dérivées en ce réel (questions (6-b-c-d))

La plupart des trinômes ont donné la représentation analytique de la fonction $f_1$ par l'usage des symboles de la valeur de la fonction $f$ et de sa dérivée en 0 par un travail s'inscrivant dans le paradigme [AC]. Par contre, deux trinômes n'ont pas essayé de déterminer cette expression. Nous pensons que l'absence de leur raisonnement revient à des contraintes du temps et notamment au manque d'expérience du travail au sein du groupe.

Par ailleurs, trois trinômes seulement ont produit un raisonnement formel complet lors de la détermination des représentations analytiques de la fonction $f_2$ par l'usage des symboles corresponds à la valeur de la fonction $f$ en 0 et de ses dérivées première et seconde en ce réel. Ils ont utilisé la méthode du répertoire didactique du Secondaire qui nécessite l'usage d'un système en vue d'identifier les termes de deux polynômes égaux. En réalité, ces étudiant ont utilisé des connaissances mobilisables du répertoire didactique du Secondaire : la méthode d'identification terme à terme (les coefficients) pour le cas d'égalité de deux polynômes de même degré. Par ailleurs, ils ont étudié cette technique Arithmetico-algébrique en deuxième année Secondaire de filière Sciences.

Nous remarquons que juste deux trinômes ont produit des pas des raisonnements de nature syntaxique afin de déterminer des représentations analytiques de la fonction $f_3$ par l'usage des symboles corresponds à la valeur de la fonction $f$ en 0 et de ses dérivées successives d'ordre 1, 2 et 3 en ce réel.

Par contre, la majorité des trinômes ont produit des raisonnements erronés lors de la détermination des approximations polynômiales locales d'ordre 2 et 3 de la fonction $f$ au voisinage de 0 par l'utilisation des symboles de sa valeur et de ses dérivées successives d'ordre 1, 2 et 3 en ce réel. Ces groupes ont effectué des *erreurs d'ordre calculatoire* lors de remplacement le coefficient du polynôme par des symboles.

A l'issue de cette analyse de notre corpus constitué des productions de 8 trinômes, nous pouvons classer les principales difficultés rencontrées par les étudiants, au point de vue de conceptions, selon trois catégories (Robert, 1998) :

- Au niveau technique : dans ce cas, ils effectuent des *erreurs d'ordre calculatoire* lors de l'usage soit la méthode d'identification terme à terme, soit la méthode de l'expression conjuguée du répertoire didactique du Secondaire dans le paradigme [AC] ;

- Au niveau des connaissances mobilisables : dans ce cas, ils effectuent des *erreurs d'ordres conceptuel et cognitif* lors de l'utilisation soit la technique du raisonnement par l'absurde dans le paradigme [AI], soit la méthode de l'étude du signe d'une fonction afin de préciser la position relative de deux courbes représentatives de deux fonctions du répertoire didactique du Supérieur dans le paradigme [AG] ;

- Au niveau des connaissances disponibles : dans ce cas, les *erreurs d'ordres conceptuel et cognitif* se réalisent lors de l'usage de la propriété de fonctions négligeables du répertoire didactique du Supérieur afin de formaliser la limite dans le paradigme [AI].

## 5.B. Analyse qualitative : exemple d'analyse des épisodes du trinôme 1

Dans cette partie, nous conduisons une analyse des différents épisodes de la séquence proposée selon les quatre axes du modèle d'analyse du raisonnement. Notre choix des épisodes repose sur la particularité des productions du trinôme 1 selon leur copie et les vidéos enregistrées. Nous pensons que l'analyse détaillée des raisonnements produits par ce groupe nous permettra d'identifier les conceptions erronées des étudiants lors de la conceptualisation des objets d'approximation locale d'une fonction. Les productions de ce trinôme sont présentées à partir de la transcription de la copie de synthèse de leurs productions, des photos prises et des vidéos enregistrées. Dans les transcriptions, on désigne par P : enseignante et E : étudiant.

Nous avons choisi de présenter l'analyse des épisodes que nous paraissent assez-important pour l'étude des conceptions (erronées ou valides) des étudiants lors de l'appropriation et l'usage raisonné du concept d'approximation locale d'une fonction.

Dans une première étape, nous présentons l'analyse selon le point de vue de l'identification du niveau de milieu et des paradigmes correspondants. Ensuite, nous identifions les registres de représentation sémiotique. Finalement, nous essayons de conduire une analyse détaillée de ces raisonnements selon du point de vue des fonctions des raisonnements, de l'usage du répertoire didactique, des analyses sémiotiques et des dimensions sémantique et/ou syntaxique.

- **Episode 1 : calculus de limite de fonctions**

Ce trinôme a produit des pas de raisonnements de nature syntaxique/sémantique lors du calculus de limite du rapport de fonctions en contrôlant la validité de la méthode utilisée. En effet, il a décomposé le quotient de fonctions sous la forme d'une somme de trois quotients. Par passage à la limite, cette transformation fait apparaitre une forme indéterminée. Alors, il mobilise la méthode de l'expression conjuguée du répertoire didactique du Secondaire.

- **Episode 2 : demande d'aide auprès de l'enseignante**

Les étudiants discutent entre eux sur la méthode à utiliser pour interpréter le résultat du calculus de limite, mais ils ne parviennent pas à mobiliser la technique adéquate pour interpréter cette expression analytique. À ce moment-là, ils demandent l'aide auprès de l'enseignante. Une étudiante de ce trinôme lève le doigt et explique à l'enseignante leur difficulté. Nous allons effectuer l'analyse de cette étape.

004 P : Est-ce que nous pouvons déterminer une relation entre la limite l d'une fonction et cette fonction ?

005 E : C'est f(x) moins l ($f(x) - l$) ?

006 P : Si tu as une fonction h tends vers l c'est-à-dire la limite de h(x) égale à l. Est-ce que h va être égale à l exactement ? Elle va être égale à quoi ? ($lim_{x \to 0} h(x) = l$)

007 E1 : L'approximation de la fonction $h$ s'annule.

008 P : Donc cette fonction $h$ admet une limite égale à l, donc h égale à quoi ?

009 E1 : Donc son approximation est proche de l.

010 P : Qu'est-ce que ça veut dire proche de l ? Comment je pourrai l'écrire ?

011E 2 : $\boldsymbol{\varepsilon}$ fois $\boldsymbol{l}$ ($\boldsymbol{\varepsilon.l}$)

012 P : Pourquoi fois $l$ ?

013 E 1 : Fois l ???

014 P : C'est $l$ lui-même

015 E1 : C'est tends vers l, c'est inférieur à l, elle va s'approcher de $l$.

016 P : Et qui t'a dit ça ? Ça peut-être plus grand que $l$.

017 E2 : Plus inférieur que $l$.

018 E : Au voisinage….

019 P : Comment on peut l'écrire au voisinage de $l$ ? Donc la fonction $h$ égale à une chose proche de $\boldsymbol{l}$.

020 E 2 : $\boldsymbol{\varepsilon}$ fois l.

021 E 1 : Pourquoi fois ?

022 P : Et $\boldsymbol{\varepsilon}$ tend vers quoi ?

023 E 2 : $(l - a)$ tend vers $l$.

024 P : Si $h$ tend vers l c'est-à-dire que la distance entre $h$ et $l$.

025 E 1 : C'est $\boldsymbol{\varepsilon}$

026 P : Et $l$, c'est très petite ?

027 E : Très petit, c'est $\boldsymbol{\varepsilon}$.

028 P : Donc la fonction $h$ qui va tendre vers $l$, c'est que la distance entre h et x est $l$ qui est en train de diminuer lorsque le x tend vers le point a, alors h(x) elle va être égale à quoi ?

029 E : …..(les étudiants ne répondent pas)

030 P : égale à $l$ ?

031 E 1 : Plus $\boldsymbol{\varepsilon}$.

032 P : Voilà. C'est plus la distance qui va s'annuler. Donc tant que la limite de la fonction $h$ égale à l. Donc on enlève le mot limite devient la fonction $h$ égale à $l$ plus la distance qui va tendre vers 0.

033 E 2 : C'est $\varepsilon$.

Tableau 4 - Transcriptions d'échanges de l'enseignante avec le trinôme 1

Ce trinôme n'a pas trouvé les moyens de faire évoluer le savoir qui est mobilisé lors de la confrontation au milieu objectif pour traduire le résultat de calcul de limite en une expression algébrique idoine reposant sur la connaissance de la notion de fonctions négligeables. En effet, au niveau de milieu $M_{-1}$, ce trinôme devrait interpréter très finement le calculus de limite et le mettre en relation avec la notion de fonctions négligeables. Il n'est pas capable de porter des réflexions sur des objets qui interviennent dans les transformations des expressions analytiques permettant le changement de paradigme. Au moment du blocage, les étudiants déclarent qu'ils ne parviennent pas à trouver la technique adéquate pour interpréter ce résultat du calcul de limite. Ils sont confrontés à une difficulté liée à la dialectique sémantique/syntaxique du concept de fonctions négligeables. Cette méthode va leur permettre de déterminer l'expression de la fonction $f$ en utilisant l'expression du calculus de limite comme un statut d'argument dans ce type du raisonnement. Dans ce cas, la propriété de fonctions négligeables est une connaissance disponible du répertoire didactique du Supérieur. À ce moment de blocage, l'enseignante intervient pour guider les étudiants de ce trinôme en vue d'utiliser la fonction $(\varepsilon(x))$ afin de formaliser la limite. Ainsi, le rôle principal de l'enseignante consiste à trouver les moyens d'intervenir pour préciser la formulation de limite. Elle place les étudiants dans une posture réflexive sur cette expression $\varepsilon(x)$ et leurs connaissances antérieures. Par son questionnement, elle tente de pallier les insuffisances du milieu. Autrement dit, il y a l'absence de rétroactions du milieu dans lequel les étudiants interagissent. Elle veut les aider à utiliser la formulation de l'expression analytique $[l + \varepsilon(x)]$, or les étudiants veulent écrire la formulation sous la forme $[l \times \varepsilon(x)]$. Cette erreur commise traduit l'existence *d'un obstacle de nature didactique* lié au choix de l'enseignante lors de l'introduction de la fonction $\varepsilon(x)$. Le rappel fait par l'enseignante, au début de la séance, sur les fonctions négligeables a un impact sur ce trinôme. En effet, il se remémore ce qui a été rappelé.

Ce trinôme est confronté à un problème résultant du changement de paradigmes de [AC] au paradigme [AI]. Il n'a pas trouvé les moyens de faire évoluer le savoir qui est entre le milieu $M_{-2}$ et le milieu $M_{-1}$. Cette *erreur d'ordre cognitif et conceptuel* est liée à la formalisation de limite par l'utilisation des connaissances disponibles "fonctions négligeables".

-   **Episode 3 : détermination de l'approximation locale d'ordre 1 de la fonction**

Lors de blocage (003E-épisode 2), les étudiants sollicitent l'enseignante. Par son intervention, elle amène les étudiants à s'interroger sur les manières de traduire ces résultats obtenus du calculus de limite. De ce fait, ils mobilisent des objets pour étudier les fonctions $f$ et $f_1$ dans le milieu de référence ($M_{-1}$). Le rôle de l'enseignante est lié aux travaux de Bloch (1999) qui a précisé le travail de professeur au niveau du milieu de référence.

Malgré cet échange, ce trinôme ne parvient pas à mobiliser la technique adéquate pour formaliser l'expression de la fonction $f$. Il n'est pas capable de porter des réflexions sur des connaissances disponibles du répertoire didactique du Supérieur permettant l'évolution de savoir entre le milieu objectif - le milieu de référence.

Nous remarquons que ce groupe est confronté au problème lié à l'utilisation de la propriété de fonctions négligeables. Après une guidance faible, il a produit un raisonnement formellement incomplet. Cet usage automatique de la méthode du registre des formules (plus précisément la propriété de fonctions négligeables) traduit l'existence d'une *erreur d'ordre cognitif et conceptuel* liée à l'appropriation du concept d'approximation locale d'une fonction. En effet, le concept de fonctions négligeables n'a de sens qu'au voisinage d'un point bien déterminé.

-   **Episode 9 : étude de la position relative de deux représentations graphiques de fonctions**

| |
|---|
| **3-c) (52)** $f(x) - f_1(x) = -\frac{1}{8}x^2 + x^2\varepsilon(x) + f_1(x) - f_1(x)$ |
| **(53)** $= -\frac{1}{8}x^2 + x^2\varepsilon(x)$ |
| **(54)** $= -\frac{1}{8}x^2 + o(x^2)$ |
| **(55)** $= \cancel{x^2(-\frac{1}{8} + \varepsilon(x))}$ |

> **(56)** Si $\varepsilon(x) - \frac{1}{8} > 0$ donc $f(x) - f_1(x) > 0$
>
> **(57)** Donc $C_f$ est au-dessus de $C_{f_1}$
>
> **(58)** Si $\varepsilon(x) - \frac{1}{8} < 0$ donc $f(x) - f_1(x) < 0$
>
> **(59)** Donc $C_f$ est au-dessous de $C_{f_1}$

Tableau 5 - Extrait de copie du trinôme 1

- <u>Identification de niveaux de milieux et de paradigmes correspondants</u>

La représentation analytique de la différence de deux fonctions prend le statut d'argument dans le milieu de référence ($M_{-1}$). Les étudiants interprètent géométriquement cette représentation analytique par l'étude de son signe. Ainsi, ils précisent les positions relatives des représentations graphiques de ces objets mathématiques : les fonctions $f$ et $f_1$. Nous pouvons dire que ce travail mathématique est piloté par les paradigmes [AI] puis [AG].

- <u>Identification de niveaux des registres de représentation sémiotique</u>

Le changement de registres est lié au paradigme [AG] à travers l'interprétation d'un calcul algébrique dans le domaine de géométrie. Ce changement induit une complexité chez ces étudiants. Ils n'ont pas la capacité d'articuler les registres analytique, algébrique, numérique et géométrique.

- <u>Du point de vue des fonctions des raisonnements, de l'usage du répertoire didactique, des analyses sémiotiques et des dimensions sémantique et/ou syntaxique</u>

**[(52)-(53)-(54)]** Décision de l'usage du signe de la différence de deux fonctions pour préciser les positions relatives de leurs représentations graphiques. Les étudiants présentent l'expression du reste sous ses deux formes : $x^2\varepsilon(x)$ et $o(x^2)$.

**[(55)]** Rétroaction sur la représentation analytique de cette différence de deux fonctions. Ce groupe donne l'expression sous forme du produit en mettant $x^2$ en facteur.

**[(56)-(57)-(58)-(59)]** Décision de l'utilisation du signe de l'expression $[\varepsilon(x) - \frac{1}{8}]$ pour déterminer les positions relatives de deux courbes.

Pour ces étudiants, la formulation de la différence de deux fonctions en produit prouve l'existence *d'erreur d'ordre cognitif et conceptuel* lors de l'utilisation des connaissances mobilisables. Ce type d'erreur est lié à l'expression du reste. Nous pouvons dire qu'il existe des difficultés récurrentes éprouvées par ce trinôme en rapport avec la conceptualisation des objets d'approximation locale d'une fonction.

En conclusion, cette analyse des raisonnements produits par des étudiants confrontés à une situation à dimension adidactique visant l'introduction de la formule de Taylor-Young, nous permet de conclure que :

- La plupart des trinômes sont confrontés aux difficultés apparaissant lors du passage du niveau du milieu objectif au milieu de référence par l'interprétation des expressions analytique et algébrique réalisées dans le milieu $M_{-2}$ ;

- Les principales difficultés des étudiants résultent du changement de paradigme soit du paradigme [AC] à celui [AI] lors de la formalisation de limite, soit du paradigme [AI] à celui [AG] lors de l'interprétation d'une représentation analytique d'une fonction ;

- Les travaux des majorités d'entre eux s'inscrivent en réalité dans le paradigme [AC] lors de la détermination des approximations polynômiales locales d'ordre 1, 2 et 3 de la fonction au voisinage de 0. Nous pouvons dire qu'il existe des conceptions liées à l'appropriation et l'usage raisonné du concept d'approximation locale d'une fonction que la plupart des étudiants ne maîtrisent pas suffisamment.

## 6. Apports et limite de la recherche

L'élaboration et la mise en œuvre d'une ingénierie didactique de développement a nécessité une collaboration étroite avec l'enseignante de la classe qui a un rôle très spécifique lors de l'enseignement du chapitre « Analyse asymptotique ». La proposition de la situation a offert aux étudiants la possibilité :

- d'interagir entre eux en vue de produire un raisonnement ;

- d'introduire eux-mêmes le concept d'approximation locale d'une fonction. La compréhension conceptuelle de la notion de développement limité déduite de la formule de Taylor-Young est obtenue en approximant une fonction aux moyens des polynômes d'une manière successive de degré 1, 2, 3, etc. ;

- d'approximer une fonction successivement par des polynômes de degré 1, 2 et 3 et de prouver son unicité au voisinage d'un réel ;

- de se positionner dans les différents paradigmes de l'Analyse standard lors de la transition Secondaire/Supérieur.

L'expérimentation aide les étudiants à percevoir qu'une fonction n-fois dérivable au voisinage d'un point peut être « approximée » par un polynôme, et sur la manière de produire l'approximation d'une fonction au voisinage d'un point et la généralisation de ce résultat dans le paradigme [AI]. Cette approximation se fait tout d'abord par une fonction polynôme de degré 1, puis par une fonction polynôme de degré 2, obtenue en ajoutant un terme à la précédente fonction, puis un travail sur la fonction polynôme de degré 3 qui est obtenue en ajoutant de nouveau terme de degré 3. Ensuite, elle explicite une manière d'établir une approximation polynômiale d'une fonction de classe $C^n$ par la détermination des valeurs de ses dérivées successives en un point donné. À partir de ces expressions analytiques et algébriques, l'enseignante a introduit la formule de Taylor-Young en s'appuyant sur les raisonnements produits par ces étudiants.

Notre travail porte la focale à la fois sur l'étude des formes de raisonnements produits par les étudiants et sur leurs discours soit au sein du groupe, soit lors de l'interaction avec l'enseignante lorsqu'ils demandent son aide. L'analyse de ces discours et de ces interactions a joué un rôle très important dans notre étude : elle nous a permis d'accéder à la nature de justification, la nature des explications, la manière dont les étudiants interagissent entre eux lorsqu'ils se posent des questions, la manière de justifier le choix des savoirs et des connaissances mobilisés lors de la production d'un raisonnement. Par ailleurs, le schéma de la structuration du milieu nous a offert la possibilité d'identifier ces savoirs, ainsi que la manière dont ils se transforment pour étudier le déroulement final de la séquence.

Dans le cas de l'élaboration de l'ingénierie didactique de développement, l'analyse des raisonnements nous a fourni des renseignements précis sur les conceptions, valides et erronées, au travers de l'utilisation que les étudiants font de leurs connaissances pour la production des raisonnements par un travail s'inscrivant dans les différents paradigmes de l'Analyse Standard ([AC], [AI] ou [AG]).

L'évolution de la construction du nouveau concept « développement limité » déduit de la formule de Taylor-Young s'effectue par les étudiants eux-mêmes à partir d'un processus mobilisant les différents registres algébrique, analytique, numérique et géométrique :

- l'expression du reste du développement limité est définie par l'interprétation du calcul de limite de fonctions en utilisant la propriété de fonctions négligeables ;

- la partie régulière du développement limité est définie par des représentations algébriques en utilisant les symboles des dérivées successives d'une fonction en un réel ;

- les approximations locales d'une fonction d'ordre 1, 2 et 3 sont déterminées via ses approximations polynômiales successives de degré 1, 2 et 3.

Par ailleurs, l'adaptation du modèle d'analyse du raisonnement, introduit par (Bloch et Gibel, 2011) et (Gibel, 2018) et développé dans notre travail de recherche, a permis d'étudier les raisonnements produits par des étudiants confrontés à une situation réelle en classe. Cette étude se réalise par la mise en considération du point de vue du fonctionnement des connaissances et des savoirs mobilisés par leurs actions sur les différents milieux, du point de vue de conceptions, du point de vue logique et du point de vue des signes. Ce modèle a permis de visualiser très précisément les connaissances et les savoirs mobilisés.

L'analyse *a priori* de la situation nous a permis de définir et de caractériser les différents niveaux de milieux, ainsi que les paradigmes correspondants. Le niveau du milieu nous donne des informations sur la nature du raisonnement produit par les étudiants confrontés à une situation à dimension adidactique. Par confrontation au niveau du milieu objectif, les étudiants réalisent un calcul syntaxique qui fait fonctionner des règles de transformations et de formulations des écritures mathématiques (des expressions analytique et algébrique). Suite à l'obtention des résultats, les étudiants doivent nécessairement trouver les moyens de les interpréter. À moment-là, ils doivent être capables de produire des réflexions sur des objets permettant de faire évoluer le savoir, notamment lorsqu'ils sont amenés à passer de la situation de référence à la situation d'apprentissage.

L'enseignante n'a pas consacré beaucoup du temps lors de cette expérimentation. Malgré ces contraintes temporelles, l'analyse de notre corpus, constitué des productions de 8 trinômes et des retranscriptions des interactions en classe, nous a permis d'identifier la nature et l'origine de leurs difficultés inhérentes à la conceptualisation des objets d'approximation locale d'une fonction en première année (PC). En effet, le modèle d'analyse des raisonnements par la prise en compte de ses composantes (fonction, dimension sémiotique, nature et justification) a joué un rôle essentiel dans l'élaboration de l'ingénierie didactique de développement, dans l'identification des conceptions des étudiants, des formes et des fonctions des raisonnements. Il a permis aussi d'identifier des obstacles liés à l'apprentissage du concept de fonctions négligeables selon deux catégories :

- Obstacle épistémologique : c'est un *obstacle d'ordre fonctionnel et topologique* qui résulte des *conceptions plus fines et plus précises* qui relèvent du paradigme [AI]. En effet, l'écriture de l'expression analytique d'approximation locale d'une fonction sans la précision de voisinage pour lequel doit être déterminée traduit l'existence des difficultés liées à la conceptualisation de cet objet mathématique. Nous pouvons dire que ces *erreurs d'ordre cognitif et conceptuel* relèvent d'un obstacle repéré dans l'histoire des mathématiques et notamment lié à la genèse historique du concept d'approximation locale d'une fonction.

- Obstacle didactique : c'est un obstacle *d'ordre culturel* lié aux choix didactiques de l'enseignement. Par exemple, l'introduction de la fonction $[x \mapsto \varepsilon(x)$, par l'enseignante, se réalise juste par le calcul de sa limite. Ce processus d'enseignement de cette fonction amène à l'existence des *erreurs d'ordre cognitif et conceptuel* aux niveaux de la formalisation de limite et l'étude du signe de l'expression de reste de l'approximation locale d'une fonction. De ce fait, les étudiants sont confrontés aux difficultés récurrentes à la mobilisation des *conceptions spécifiques* liées au changement de paradigmes. Il existe une rupture conceptuelle qui s'opère entre le niveau Secondaire et le niveau Supérieur.

Nous pouvons dire que certains obstacles de nature épistémologique, didactique et culturelle résultent du changement de paradigme lors de la transition Secondaire/Supérieur.

## 7. Conclusion et perspectives

Dans ce travail de recherche, nous nous intéressons, dans le cadre de la théorie des situations didactiques (Brousseau, 1998), à la construction des raisonnements produits par les étudiants, ainsi qu'à la question du contrôle de ces raisonnements afin d'identifier la nature et l'origine des erreurs effectuées par les étudiants lors de la conceptualisation de la notion d'approximation locale d'une fonction. Pour ceci, nous avons construit et mis en œuvre l'ingénierie didactique de développement qui a nécessité une collaboration étroite avec l'enseignant. Elle a un rôle très spécifique en classe. Elle s'est investie dans l'élaboration du chapitre « Analyse Asymptotique » dans lequel la situation à dimension adidactique offre aux étudiants la possibilité de faire fonctionner leurs connaissances et leurs savoirs du répertoire didactique par leurs actions successives sur les différents milieux afin d'accéder au sens du concept d'approximation locale d'une fonction dans le paradigme [*Analyse Infinitésimale*]. L'enseignante parvient, par une guidance faible (Bartolini-Bussi, 2009), à aider les étudiants de produire des raisonnements originaux permettant de trouver leur place dans la construction du concept d'approximation locale d'une fonction. En fait, dans la situation, ils accèdent à ce concept mathématique par une succession de questions qui les amènent à élaborer des raisonnements à partir du concept de limite afin de faire le lien avec les notions de fonctions négligeables et arriver ensuite aux différentes écritures syntaxiques des approximations polynômiales d'une fonction au voisinage d'un réel par la mobilisation des registre analytique, algébrique, numérique et géométrique.

Nous avons conduit une analyse *a priori* et une analyse *a posteriori* de l'expérimentation dans le cadre de la théorie des situations didactiques afin de produire des éléments de réponses à nos questions de recherche et à elles proposées par l'enseignante sur la nature et l'origine des difficultés rencontrées par les étudiants lors de la conceptualisation des notions d'approximation locale d'une fonction en première année (PC). Notre formulation des questions en tant que didacticien diffère de celle proposée par l'enseignante. Nous nous appuyons sur des concepts didactiques, et notamment la notion d'obstacle, (Brousseau, 1998).

Nous pouvons dire que l'ingénierie didactique de développement a permis de pratiquer le raisonnement afin de donner du sens à la conception du concept d'approximation locale d'une fonction dans le paradigme [*Analyse Infinitésimale*] et de percevoir ses différentes fonctions dans l'apprentissage : décider de l'utilisation d'une connaissance (au niveau technique, connaissances mobilisables ou disponibles), expliquer, contrôler, justifier, convaincre, réfuter, prouver, démontrer, etc.

L'expérience d'apprentissage, inhérente à la situation, a été la première rencontre des étudiants avec la notion d'approximation d'une fonction à partir du concept « développement limité d'une fonction au voisinage d'un réel ». Cette expérimentation a atteint son but : offrir la possibilité aux étudiants de produire un raisonnement de nature sémantique pour prouver l'unicité de l'approximation polynômiale d'ordre 2 d'une fonction au voisinage d'un réel et un raisonnement articulant les dimensions sémantique et syntaxique par des interprétations du calculus de limite afin d'introduire la formule de Taylor-Young. Cette ingénierie contribue à favoriser une approche analytique permettant l'introduction de la définition formelle du concept de développement limité dans le paradigme [AI].

Dans notre travail, la mise en œuvre d'une situation à dimension adidactique, couplée au modèle d'analyse des raisonnements, est un dispositif intéressant qui nous a permis de repérer les conceptions (valides et erronées) des étudiants au travers de leurs productions.

Cette recherche a conduit à accorder une place importante à la conceptualisation de la notion de « limite » d'une fonction, au niveau Supérieur, qui est un élément clé afin d'aider les étudiants à accéder au sens du concept d'approximation locale d'une fonction dans le paradigme [AI] en première année (PC). Cette étude pourra jouer un rôle important dans la formation des enseignants au point de vue de la construction des nouvelles connaissances au début du cursus dans le Supérieur.

# Références

Fatma Belhaj Amor
Laboratoire de mathématiques et de leurs applications (LMAP)
*e-mail:* f.belhaj-amor@univ-pau.fr